\documentclass[11pt,a4paper,reqno]{amsart}
 \usepackage{amsfonts,amsthm, amscd, epsfig, amsmath, amssymb,enumerate}
\usepackage[dvips]{psfrag}
\numberwithin{equation}{section}%

\makeatletter
\def\captionfont@{\footnotesize}
\def\captionheadfont@{\scshape}

\long\def\@makecaption#1#2{%
  \vspace{2mm}
  \setbox\@tempboxa\vbox{\color@setgroup
    \advance\hsize-6pc\noindent
    \captionfont@\captionheadfont@#1\@xp\@ifnotempty\@xp
        {\@cdr#2\@nil}{.\captionfont@\upshape\enspace#2}%
    \unskip\kern-6pc\par
    \global\setbox\@ne\lastbox\color@endgroup}%
  \ifhbox\@ne 
    \setbox\@ne\hbox{\unhbox\@ne\unskip\unskip\unpenalty\unkern}%
  \fi
  \ifdim\wd\@tempboxa=\z@ 
    \setbox\@ne\hbox to\columnwidth{\hss\kern-6pc\box\@ne\hss}%
  \else 
    \setbox\@ne\vbox{\unvbox\@tempboxa\parskip\z@skip
        \noindent\unhbox\@ne\advance\hsize-6pc\par}%
\fi
  \ifnum\@tempcnta<64 
    \addvspace\abovecaptionskip
    \moveright 3pc\box\@ne
  \else 
    \moveright 3pc\box\@ne
    \nobreak
    \vskip\belowcaptionskip
  \fi
\relax } \makeatother
\def\writefig#1 #2 #3 {\rlap{\kern #1 truecm
\raise #2 truecm \hbox{#3}}}

\newtheorem{theorem}{Theorem}[section]
\newtheorem{lemma}[theorem]{Lemma}
\newtheorem{prop}[theorem]{Proposition}
\newtheorem{cor}[theorem]{Corollary}

\newcommand{\LL}{{\mathbb L}}

\newcommand{\CC}{{\mathbb C}}
\newcommand{\RR}{{\mathbb R}}

\newcommand{\NN}{{\mathbb N}}
\newcommand{\PP}{{\mathbb P}}

\newcommand{\EE}{{\mathbb E}}

\newcommand{\ep}{{\epsilon}}

\newcommand{\Av}{{\rm Av}}  

\newcommand{\II}{{\mathbb I}}

\newcommand{\dto}{\downarrow}
\newcommand{\uto}{\uparrow}

\newcommand{\cA}{\ensuremath{\mathcal A}}
\newcommand{\cB}{\ensuremath{\mathcal B}}

\newcommand{\cE}{\ensuremath{\mathcal E}}

\newcommand{\cG}{\ensuremath{\mathcal G}}

\newcommand{\cM}{\ensuremath{\mathcal M}}

\newcommand{\cP}{\ensuremath{\mathcal P}}

\newcommand{\cS}{\ensuremath{\mathcal S}}

\let\a=\alpha \let\b=\beta   \let\d=\delta  
 \let\g=\gamma     \let\k=\kappa  \let\l=\lambda
      \let\o=\omega    \let\p=\pi  
  \let\s=\sigma \let\t=\tau   
  
   \let\G=\Gamma  \let\L=\Lambda 
\let\O=\Omega \let\P=\Pi     

\let\Z=\integer

\let\neper=e
\let\ii=i

\let\sset=\subset

\let\<=\langle
\let\>=\rangle

\outer\def\nproclaim#1 [#2]#3. #4\par{\medbreak \noindent
   \talato(#2){\bf #1 \Thm[#2]#3.\enspace }%
   {\sl #4\par }\ifdim \lastskip <\medskipamount
   \removelastskip \penalty 55\medskip \fi}
\def\thmm[#1]{#1}
\def\teo[#1]{#1}
\def\sttilde#1{%
\dimen2=\fontdimen5\textfont0
\setbox0=\hbox{$\mathchar"7E$}
\setbox1=\hbox{$\scriptstyle #1$}
\dimen0=\wd0
\dimen1=\wd1
\advance\dimen1 by -\dimen0
\divide\dimen1 by 2
\vbox{\offinterlineskip%
   \moveright\dimen1 \box0 \kern - \dimen2\box1}
}

\def\\{\hfill\break}

\def\tthsp{\kern .083333 em}

\def\?{\mskip -10mu}

\begin{document}

\title[Spectral characterisation of ageing]
{Spectral characterisation of ageing: the REM-like trap model}
\author{Anton Bovier}
\address{Weierstrass Institut f\"ur Angewandte Analysis und
  Stochastic, Mohrenstrasse 39, 10117 Berlin, Germany, and
Mathematisches Institut, Technische Universit\"at Berlin, Strasse des
17. Juni 136, 10623 Berlin, Germany. email: bovier@wias-berlin.de}
\author {Alessandra Faggionato}
\address{Weierstrass Institut f\"ur Angewandte Analysis und
  Stochastic, Mohrenstrasse 39, 10117 Berlin, Germany. email:
  faggiona@wias-berlin.de}
\begin{abstract}{We review the ageing phenomenon in the context of simplest
  trap model, Bouchaud's REM-like trap model from a spectral theoretic
  point of view. We show that the generator of the dynamics of this
  model can be diagonalised exactly. Using this result, we derive
  closed expressions for correlation functions in terms of complex
  contour integrals that permit an easy investigation into their large
  time asymptotics in the thermodynamic limit. We also give a `grand
  canonical' representation of the model in terms of the Markov
  process on a Poisson point process . In this context we analyse the
  dynamics on various time scales.  }\\

\end{abstract}
\thanks{ Research supported in part by the DFG in the 
Dutch-German Bilateral Research Group  "Mathematics of Random Spatial 
Models from Physics and Biology".}

\maketitle
\section{Introduction}

The particular properties of the long term dynamics of 
many complex and/or disordered systems have been the subject of great
interest in the physics, and, increasingly, the mathematics
community. The key paradigm here is the notion of \emph{ageing}, a
notion that can be characterised in terms of scaling properties of
suitable autocorrelation functions. Typically, ageing can be
associated to the existence of \emph{infinitely many} time-scales that
are inherently relevant to the system. In that respect ageing systems
are distinct from \emph{metastable} systems which are characterised by
a finite number of well separated time scales corresponding to the
live times of different metastable states. 

Ageing systems are rather difficult to analyse, both numerically and
analytically.  Most analytical results, even on the heuristic level,
concern either the Langevin dynamics of spherical mean field spin
glasses, or \emph{trap models}, a class of artificial Markov processes
that in some way is to mimic the  long term dynamics of highly
disordered systems
(see e.g. \cite{BCKM}). 

One of the natural questions one is led to ask when being confronted
with phenomena related to multiple time scales is whether and how they
can be related to \emph{spectral properties}. This relationship has
been widely investigated in the context of Markov processes with
metastable behaviour (see e.g. \cite{D1,D2,D3,FW,GS,BEGK}), and it would
be rather interesting to obtain a spectral characterisation of ageing
systems as well, at least in the context of Markov processes. 
To our knowledge, this problem has not been widely studied so far. The
only paper dealing with the problem is the paper \cite{BM} by Butaud
and Melin that have tackled one of the simplest trap models and on
which we will comment below, and  
\cite{FIKP} and \cite{M} that investigate convergence to equilibrium
in the Random Energy Model (REM).

The present paper is intended to make a modest first step into this
 direction by analysing the relation between spectral properties
and ageing  rigorously in the REM-like trap model. While this model
may seem misleadingly simple, it has in the past provided valuable
insights into the mechanisms of ageing, and it is our hope that the
analysis presented here will provide useful guidelines for further
investigations of more complicated models.

The paper will be divided into two parts. In the first we analyse the 
REM-like trap model in the standard formulation of Bouchaud
\cite{Bou}. In the second part we go a step further and reformulate
the model in a slightly different way as a Markov process on a Poisson
point process. This formulation makes the relation to the real REM
more suggestive (see \cite{BBG1, BBG2} for a full analysis), and
allows in a 
natural way to study the dynamics of the model on different time
scales. 

\section{The REM--like trap model}\label{inizio}

Let us begin to recall the definition of \emph{trap models} as
introduced by Bouchaud and Dean\cite{Bou}. 
 Let $\cG=(\cS, \cE) $ be a finite graph:  $\cS$
denotes the set of vertices   and  $\cE$ denotes the  set of edges,
$\underline E:=\{E_i,\, i\in \cS\}$ be a random field,
 called \emph{energy landscape} and let   $Y(t)$ be a
continuous--time random walk  on $\cG$  with $\underline
E$--depending transition rates $c_{i,j}$: i.e. $c_{i,j}>0$ iff
$\{i,j\}\in\cE$ and
$$ \PP( Y(t+dt)=j\,|\, Y(t)=i)= c_{i,j} dt.$$
Setting    $\t_i:=\sum _{j\not = i} c_{i,j} $ and $p_{i,j}:=
c_{i,j} /\t_i$, the random walk $Y(t)$ can be described as
follows: after reaching the site $i$ the system waits an
exponential time with expectation $\t_i$ and then it jumps to an
adjacent site  $j$  with probability $p_{i,j}$. In the trap model,
the transition rates are assumed to satisfy the following
properties:
\begin{align}
&e^{E_i}c_{i,j}=e^{E_j}c_{j,i},\quad \forall  \{i,j\}\in\cE,\\
&\EE(\t_i)=\infty
\end{align}
where  $\EE$ denotes the expectation w.r.t. the random field 
$\underline E$. Since in several physical experiments (see
\cite{bibo}) the system initially in  equilibrium at high temperature 
$T\gg T_g$  is quickly cooled under the transition temperature 
$T_g$ and then its response to an external perturbation is measured,
 it is reasonable to consider $Y(t)$ with uniform
initial distribution.\\
A classical time--time correlation function is given by
$$
\Pi(t,t_w):=\PP(\,Y(s)=Y(t_w), \quad \forall s\in [t_w,t_w+t])
$$
In order to observe ageing  it is necessary to consider a thermodynamic
limit, with the size of $\cG$ going to infinity, and possible a
suitable  time--rescaling. Rather recently, there have been a number
of rigorous papers devoted to the analysis of trap models on the
lattices $ \Z$ \cite{FIN,FINS2,BC} and $\Z^d$ \cite{BCM,C}.

In this paper we consider the simples trap model, called the 
\emph{REM--like trap model} \cite{Bou} that corresponds to choosing
$\cG$ to be the complete graph on $N$ vertices
 i.e.
$$ 
\cG_N=(\cS_N, \cE_N),\;\;\; \cS_N:=\{1,2,\dots, N\},\;\;\;
\cE_N:= \{\,\{i,j\}\; i\not =j\in \cS_N\},
$$
 and to take as energy 
landscape  a family  $\underline E=\{E_i\,:\,i\in\NN\}$ of 
 independent exponential  random variables with parameter $\a$ such that 
 $0<\a<1$.  
Given $N\in \NN$, let $Y_N(t)$ be the continuous--time  random walk on
$\cG_N$  with  transition
rates $c_{i,j}= e^{-E_i}/N$ for $i\not = j$. 
 Setting $x_i= e^{-E_i}$ the infinitesimal generator of
 the random walk is given by
\begin{equation}\label{perleaiporci}
\LL_N :=
\begin{pmatrix}
\frac{(N-1)x_1}{N} & -\frac{x_1}{N} & \ldots &-  \frac{x_1}{N} \\
-\frac{x_2}{N} & \frac{(N-1) x_2}{N} & \ldots & - \frac{x_2}{N} \\
\vdots        & \vdots        & \ddots & \vdots \\
-\frac{x_N}{N} &- \frac{x_N}{N} & \ldots &  \frac{(N-1) x_N}{N}
\end{pmatrix}
\end{equation}
The dynamics can be described as follows: after reaching the state $i$
the system waits an exponential time of mean $\frac{N}{N-1}e^{E_i} $ and then
jumps with uniform probability to another state.  Although strictly
speaking the mean waiting time is given by $\frac{N}{N-1}e^{E_i} $, we
call $\t_i:=x_i^{-1}=e^{E_i}$ waiting time (the discrepancy is
negligible in the thermodynamic limit $N\uto\infty$).\\
Note that  $\t_i$ and $x_i$ have distributions respectively given
by
$$
p(\t)d\t=
\a\t^{-1-\a}d\t \quad (\t\geq 1);\qquad
p(x)dx  =\a x^{\a-1}d x \quad  ( 0<x\leq 1 ),
$$
in particular $\EE(\t_i)=\infty$. Moreover, the equilibrium measure is given by
$\mu_{eq}(i)=\t_i/ \bigl(\sum_{j=1}^N \t_j\bigr)$. We are interested
in the out--of--equilibrium dynamic with  uniform initial
distribution. $\PP_N$ denotes   the law of this random walk given a
realization of the random variables $E_i$.

\bigskip

Ageing in the REM-like trap model is manifest from the asymptotic
behaviour 
of the time--time correlation function 
\begin{equation}\label{corfun}
\Pi_N (t,t_w):=
\PP_N\bigl(\,Y_N(s)=Y_N(t_w),\; \forall s\in [t_w, t_w+t]\bigr)
\end{equation}
Namely, as shown in \cite{Bou},
 for almost all $\underline E$ and for all $\theta>0$
\begin{equation}
\lim_{t_w \uto\infty} \lim_{N\uto\infty}  \Pi_N (\theta t_w ,t_w) =
 \frac{\sin (\pi \a)}{\pi}  \int_{\frac \theta {1+\theta }}^1 u^{-\a}(1-u)^{\a-1}du
\end{equation}

 Our main aim  here is to show that  the ageing behaviour of the system, 
derived  in \cite{Bou} using renewal arguments,  can be obtained 
solely from \emph{spectral information} about the
generator $\LL_N$.
The method  developed below will allow us to get further information
on $Y_N(t)$ from the spectral properties of $\LL_N$. In particular,
given a function $h$ on $(0,\infty)$, it is possible to describe the
asymptotic behaviour of $\EE_N (h (x_N(t))$ and
$\EE_N (h (\t_N(t))$, where $\EE_N$ denotes the expectation w.r.t. 
$\PP_N$ and  $x_N(t), \t_N(t)$ are defined  as
$$ x_N(t)=x_k,\;\; \t_N(t)=\t_k\;\; \text{ if }\;\; Y_N(t)=k .$$
These results will allow us to investigate the property that the
system with high probability  visits deeper and deeper traps. i.e.
sites with larger and larger
 waiting time $\t_i$ (see subsection \ref{step2}).

We start by giving a complete description of the  eigenvalues and
 eigenvectors of $\LL_N$. Let $\mu=\mu_N$
 be the measure on $\cS_N $ with
$\mu(i)= x_i^{-1}= \t_i$. Note that  $\LL_N$ is a symmetric
operator on $L^2(\mu)$ and trivially  $\LL_N \II=0$  where $\II$
is the vector  with all entries equal to $1$. 
The following proposition is based on elementary linear algebra:
\begin{prop}\label{structure1}
Let   $x_1, x_2, \dots ,x_N$  be all distinct. Then, 
$\LL_N$ has $N$ positive  simple  eigenvalues $0=\l_1<\l_2\dots<\l_N$
such that
$$
\{ \l_1,\l_2,\dots,\l_N\}=\{\l\in\CC\,:\, \phi(\l)=0 \}
$$
where $\phi(\l)$ is the  meromorphic  function
\begin{equation}\label{chiave}
\phi(\l):=\sum_{j=1} ^N \frac{\l}{x_j -\l},\qquad (\l\in \CC).
\end{equation}
If  the $x_i$ are labelled such that
 $x_1<x_2<\dots < x_N$, then $ x_i<\l_{i+1} < x_{i+1}$ for
$i=2,\dots, N$. Moreover, for any $i=1,\dots, N$,  the vector
$\psi^{(i)}\in \RR^N$ defined as
$$
\psi^{(i)} _j:=\frac{x_j}{x_j-\l_i},\qquad \text{for } j=1,\dots,N
$$
is an    eigenvector of $\LL_N$ with eigenvalue $\l_i$.
$\psi^{(1)}, \dots, \psi^{(N)}$ form  an orthogonal basis of
$L^2(\mu)$.
\end{prop}

Since the $x_i$ have a absolutely continuous distribution, we
trivially have the 
\begin{cor}\label{trivialcor}
The assertions of Proposition \ref{structure1} hold with probability
one for all $N$.
\end{cor}

\begin{proof}
Let $\l$ be a generic eigenvalue and let us write a related
eigenvector 
$\psi$ as $\psi= a (1,\dots, 1)^t+w$ where $\sum _{j=1}^N w_j =0$. Since
$(\LL_N \psi )_j= x_j w_j$, we have to solve the system
\begin{equation}\label{tata}
x_j w_j = \l a +\l w_j, \qquad \forall j=1,\dots, N.
\end{equation}
Since $x_1,\dots, x_N$ are distinct, 
it must be true that $a\not=0$ (otherwise we get $\psi=0$). Without
loss of 
generality, we set $a=1$.  Note that $\l\not = x_j$ for $j=1\dots N$, 
as otherwise (\ref{tata}) would imply $\l=0=x_j$. Therefore we get
 $w_j= \frac{\l}{x_j-\l} $. Since it must be true that 
 $\sum _{j=1}^N w_j=0$, we
get that $\l$ is an eigenvalue and  $\psi$, with $\psi_j
=\frac{x_j}{x_j-\l}$, 
the  
corresponding  eigenvector, iff $\phi(\l)=0$. This imply that $\phi$ has at
most $N$ zeros.
Since $ \phi(0)=0$, and, for real $\l$, $\lim _{\l \dto x_i}\phi(\l)=-
\infty$, $\lim_{\l\uto x_i}\phi(\l)=\infty$,  we get that $\phi$ has
exactly $N$ zeros. The conclusion of the proof is trivial.
\end{proof}

Proposition \ref{structure1} has the following simple
corollary.

\begin{cor}\label{cor.1}
With probability one,
 the spectral distribution  $\s_N:=\Av _{j=1}^N \d_{\l_j}$ converges 
weakly to the measure $\a x^{\a-1} dx$ on $[0,1]$.
\end{cor}

\noindent {\bf Remark.} The results of Proposition \ref{structure1}
are incompatible with the predictions of \cite{BM}. The discrepancy is
particularly pronounced in the case of the eigenfunction. The reason
for this is an inappropriate use of perturbation expansion in
\cite{BM}. We will explain this in some detail in an appendix. 

We will now show that Proposition \ref{structure1} allows to derive
the asymptotics of the autocorrelation functions easily. In fact, it
contains far more information on the long time behaviour of the
systems some of which we will bring to light later.

Recall that  $p_t(i,j)$, the probability  to jump from  $i$ to $j$ in 
an interval of  time $t$, can be expressed as
$p_t(i,j)= \bigl (e^{-t\LL_N}\bigr)_{i,j}.$ In particular, by writing 
$\nu_t$ for the  probability distribution  of $Y_N(t)$  and thinking of
the  Radon derivative $\frac{d\nu_t}{d\mu}$    as column vector,
$$
\frac{d\nu_t}{d\mu}= e^{-t\LL_N} \frac{d\nu_0}{d\mu},
$$
thus  implying
\begin{equation}\label{dariofo}
\frac{d\nu_t}{d\mu}=\sum _{k=1}^N  \frac{<  \frac{d\nu_0}{d\mu} ,
\psi^{(k)}>}{<\psi^{(k)},\psi^{(k)}>}  e^{-t\l_k}\psi^{(k) }
\end{equation}
The above formulas are true for an arbitrary initial distribution. 
Taking $\nu_0$ the uniform distribution,  by Proposition \ref{structure1} we get
$$ 
\frac{d\nu_o}{d\mu}=\sum _{k=1}^N \g_k \psi^{(k)}, \text{ where }
 \g_k^{-1} :=
<\psi^{(k)},\psi^{(k)}>= \sum_{j=1}^N  \frac{x_j}{(x_j-\l_k)^2}.
$$
Then, by  Proposition \ref{structure1} and  (\ref{dariofo}),
\begin{align}
& \Pi_N (t,t_w)=
\sum_{j=1}^N\sum_{k=1}^N \frac{\g_k e^{-\l_k t_w}}{x_j-\l_k}
e^{-\frac{N-1}{N}x_j t} \label{int1}\\
& \EE_N\bigl(h(x_N(t))\bigr)=
\sum_{j=1}^N  \sum_{k=1}^N
 \frac{\g_k e^{-\l_k t}}{x_j-\l_k} h(x_j) \label{cap1}
\end{align}

The above formulas (that may appear rather ugly at first sight) admit
a nice complex
integral representation  through  the following lemma.
\begin{lemma}\label{rule}
 Let  $\g$ be a positive oriented loop  on $\CC$  containing in its interior
$\l_1, \dots, \l_N $. Let  $g$ be  an holomorphic function on a domain
$D\sset \CC$ with $\g\sset D$. 
  Then, for any $j=1\dots, N$,
\begin{equation}
 \sum_{k=1}^N \frac{\g_k  g(\l_k)}{x_j-\l_k} =\frac{1}{2\pi i}
\int_\g \frac{g(\l)}{\phi(\l)(x_j-\l) } d\l.
\end{equation}
\end{lemma}
\begin{proof}
Let us set  $X:=\{x_1,\dots,x_N\}$ and  $\L:=\{\l_1,\l_2,\dots,
\l_N\}$. 
Then $\phi(\l) $ is an holomorphic function on $\CC\setminus X$, where
$\phi'(\l)= \sum_{j=1}^N \frac{x_j}{(x_j-\l)^2}$, and in particular 
$\phi'(\l_j)= \g_j^{-1}$. Moreover, the function $[\phi(\l)
(x_j-\l)]^{-1}$ that is  a priori defined on $\CC\setminus (X\cup\L)$ 
 can be analytically  continued to $X$ to a meromorphic function with
 simple  poles only  at the
points of $\L$. Now the conclusion follows from a
trivial application of the residue theorem. 
\end{proof}

We can obviously use Lemma \ref{rule} to rewrite   formulas
(\ref{int1}) and (\ref{cap1}) in the form
\begin{align}
& \Pi_N (t,t_w)=
 \frac{1}{2\pi i}\int_\g\frac{e^{- t_w \l }}{\l}
\Bigl( \Av_j\frac{e^{-\frac{N-1}{N}x_j t}}{x_j -\l} \Big{/}
\Av_j \frac{1}{x_j-\l}\Bigr)
 d\l \label{integrale}\\
& \EE_N \bigl(h(x_t)\bigr)=
\frac{1}{2\pi i}\int_\g\frac{e^{-t\l}}{\l}
\Bigl(\Av_j\frac{h(x_j)}{x_j-\l} \Big{/}\Av_j \frac{1}{x_j-\l}\Bigr)d\l
\label{caporale}
\end{align}
where $\Av _j $ denotes the average over $j=1,2,\dots, N$.

The above integral representations of $\Pi_N (t,t_w)$  and
$ \EE_N \bigl(h(x_t)\bigr)$  have two main advantages.
First, the appearance of averages allows to compute   their limiting behaviour as
$N\uto\infty$  easily  by using the ergodicity of
the random field $\underline E$. Second, by means of the residue
theorem, their Laplace transform  can be easily computed in order
to derive  the asymptotic behaviour of $\Pi_N (t,t_w)$ and  
$ \EE_N \bigl(h(x_t)\bigr)$  for $N,t_w,t\gg 1$ 
(see subsections \ref{step1}, \ref{step2}).

 A much more general derivation of the above integral representations is discussed  in  Appendix \ref{rinnovo}.

\subsection{Ageing behaviour of $\Pi_N(t, t_w)$}\label{step1}
\begin{prop}\label{schneider}
Let us define
\begin{equation}\label{iacopo11}
\P (t,t_w):=\frac{1}{2\pi i}\int_\g \frac{e^{- t_w \l}}{\l}
\frac{ \EE_x \bigl (\frac{e^{-xt}}{\l-x} \bigr)}{\EE_x\bigl
 (\frac{1}{\l-x}  \bigr)} d\l.
\end{equation}
where $\EE_x$ is the expectation w.r.t. the measure $\a x^{\a-1}dx$ on
$[0,1]$ and $\g$ is any positive oriented complex 
loop  around the interval $[0,1]$. Then
\begin{equation}
\lim _{N\uto\infty} \P_N (t,t_w)= \P (t,t_w)
 \qquad \forall t,t_w,\qquad a.s.\,.
\end{equation}
\end{prop}
\begin{proof}
Recall (\ref{integrale}) and fix $0<\d<1/2$. Due to analyticity,
we can choose the integration contour  $\g$ to
have  distance $1$ from the
segment $[0,1]$. For  each $\l\in\g$, the random variables
$(x_j-\l)^{-1}$, $j\in\NN$,  are i.d.d. and bounded. Therefore,
for a suitable positive constant $c>0$,
\begin{equation}\label{muffin}
\PP\Bigl(\Bigl |
 \Av_{j=1}^N \frac{1}{x_j-\l} -\EE_x\bigl
 (\frac{1}{\l-x}  \bigr) \Bigr | \geq N^{-\frac{1}{2}+\d} \bigr )\leq
e^{ -c\, N^{2\d} } \qquad \forall \l \in \g.
\end{equation}
Since  for each   $x\in [0,1]$ and $\l\in\g$,  $|
\frac{\partial}{\partial \l} (x-\l)^{-1} |\leq 1$, a simple chaining
argument  allows  to
deduce from the pointwise estimate (\ref{muffin}) uniform control in
$\l$. With the   Borel--Cantelli lemma  one can then infer that,
a.s.,
\begin{equation}
\sup _{\l\in\g} \,\Bigl |
 \Av_{j=1}^N \frac{1}{x_j-\l} -\EE_x\bigl
 (\frac{1}{\l-x}  \bigr) \Bigr | \leq c\,N^{-\frac{1}{2}+\d}, \qquad  
\forall N\in \NN.
\end{equation}
Similar arguments show  that, a.s., given $M\in \NN$ there exists
 a constant $c_M$ such that  
\begin{equation}
\sup_{M-1\leq t\leq M} \sup _{\l\in\g}\, \Bigl |
 \Av_{j=1}^N \frac{e^{-\frac{N-1}{N} x_j t} }{x_j-\l} -\EE_x\bigl
 (\frac{e^{-x_j t} }{\l-x_j }  \bigr) \Bigr | \leq c_M\,
 N^{-\frac{1}{2}+\d},
 \qquad   \forall N \in \NN.
\end{equation}
Note that, for each $\l\in \g$,  $\Av_{j=1}^N (x_j-\l)^{-1}$ is a
 convex combination of points of modulus larger or equal than $1/2$, 
contained in a angular sector with angle non larger than a suitable constant  $ c<\pi$.
 In particular, $\bigl |\Av_{j=1}^N (x_j-\l)^{-1}\bigr|\geq  c'>0$  for all
$N$.
 From here the assertion of the proposition 
 follows from Lebesgue's Dominated Convergence Theorem.
\end{proof}

Given $\theta>0$ we are interested in the limit of $\Pi  (\theta t_w, t_w)$ as
$t_w\uto \infty$. This will be done using the Laplace transform
of $ \Pi  (\theta t_w, t_w)$,
\begin{equation*}
\hat \Pi(\theta,\o):= \int_0^\infty  e^{-\o t_w}\Pi(\theta t_w,t_w)   
d t_w,\qquad (\Re(\o)>0 ).
\end{equation*}
The computation of this Laplace transform is trivial if we use 
 the integral expression (\ref{iacopo11}).

\begin{figure}[!ht]
    \begin{center}
      \psfrag{a}[l][l]{$-\Re(\o)$}
      \psfrag{b}[l][l]{$0$}
     \psfrag{c}[l][l]{$1$}
     \psfrag{d}[l][l]{$\gamma$}
           \includegraphics[width=7cm]{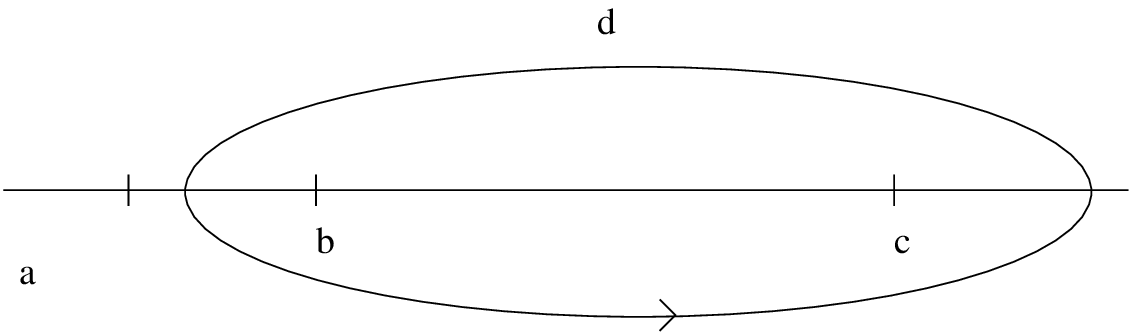}
    \label{patti}
    \end{center}
  \end{figure}

Let $\o\in \CC$ with $\Re (\o)>0$ and  fix a
positive oriented loop $\g$ around the segment $[0,1]$,  such that
$\g\sset \{ z\in \CC: \Re(z)>-\Re(\o)\}$ as in fig. \ref{patti}.
Then, $\Re(\o+\l+x\theta)>0$, for $x\in[0,1]$ and $\l\in\g$,
so that (\ref{iacopo11}) implies
$$
\hat \Pi(\theta,\o)=
\EE_x \Bigl (
 \frac{1}{2\pi i}\int_\g \Bigl [
\l (\l-x) (\l+\o+\theta x) \EE_{\bar x }(\frac{1}{\l-\bar x}) \Bigr ]^{-1}d\l \Bigr ).
$$
Here $\EE_x$ and $\EE_{\bar x}$ denote the expectation w.r.t. to the
measure $\a x^{\a-1}dx $ on $[0,1]$.

Let us consider the change of variable   $ z=\frac{1}{\l}$ and  write
$\hat \g$ for the path $\g$ with inverted orientation (i.e., positive
oriented w.r.t. $\l=\infty$).  Then we get
$$
\hat \Pi(\theta,\o)=
\EE_x \Bigl (
 \frac{1}{2\pi i}\int_{\hat \g} \Bigl [
 (1-zx) (1+z\o+z\theta x) \EE_{\bar x }(\frac{1}{1-z\bar x}) \Bigr ]^{-1} dz\Bigr ).
$$
Given $x\in [0,1]$,  the integrand is a meromorphic function  in 
$\CC \setminus [1, \infty)$ that has only a single  pole of order $1$ 
    inside $\hat \g$, namely at
$z=-(\o+x\theta)^{-1}$. By the  residue
theorem we get
\begin{equation}\label{davide}
\hat \Pi(\theta,\o)
 =\EE_x \Bigl ( \frac{1}{\o+x\theta +x}\Big{/} \EE_{\bar x}\bigl(
\frac{\o+x\theta}{\o+x\theta+\bar x} \bigr ) \Bigr).
\end{equation}

\begin{lemma}\label{pisa}
The r.h.s. of (\ref{davide}) is well defined and holomorphic for any
$\o \in \CC\setminus (-\infty, 0]$. In particular, the function 
$\hat \Pi (\theta,\o)$, defined for $\Re (\o)>0$, can be
analytically continued to the set $\CC\setminus (-\infty,0]$.
\end{lemma}
\begin{proof}
As proved in \cite{doe1}, Chapter 3,
 the Laplace transform $\hat \Pi (\theta,\o)$ is holomorphic 
on the set of convergence points. Therefore, we only need to 
show that the r.h.s.  of (\ref{davide}) is well defined and
 holomorphic on $\Im(\o)\not=0$. Let us assume   $\Im(\o)>a>0$. 
Then, trivially, $\forall x, \bar x\in [0,1]$,
$$
\frac{\o+x\theta}{\o+x\theta+\bar x}\in \cB:=
\{z\in\CC\,:\, z=|z| e^{i\theta}\text{ with } 0\leq 
\theta \leq \theta_0 ,\;\; |z|\geq c \}
$$
for  suitable constants $c,\theta_0$ depending on $a$ and such that
$\theta_0 <\pi$.
Moreover,
since $\lim _{|\o|\uto\infty} \frac{\o+x\theta}{\o+x\theta+\bar x}=1$,
\begin{equation}\label{palline}
0<c_1(a)\leq \bigl |\frac{\o+x\theta}{\o+x\theta+\bar x}
\bigr|\leq c_2(a) \qquad
\forall a>0, \forall \o :\, \Im(\o)\geq a
\end{equation}
By  (\ref{palline}) and the geometry of  $\cB$, we have that
$\EE_{\bar x}\bigl(\frac{\o+x\theta}{\o+x\theta+\bar x} \bigr )$
is well defined and has distance
$c_3(a)>0$ from the origin.  Moreover,
$|\o+x\theta+x  |\geq \Im(\o)$. Therefore,  the r.h.s.
  of (\ref{davide}) is well defined and, 
due to the previous estimates and  Lebesgue's dominated convergence theorem, 
it is continuous on $\{\Im(\o)\not =0\}$, thus implying continuity on 
$\CC\setminus (-\infty,0]$.

We recall Morera's theorem: if $f(\o)$ is defined and continuous in a
open set 
$\O\sset \CC$ and if $\int_\g f \,d\o=0 $ for all closed curves 
$\g$ in $\Omega$, then  $f(\o)$ is holomorphic in $\O$.
Therefore, using   Fubini's  and Morera's theorems, 
 one can prove that  the function 
 $\EE_{\bar x}\bigl( \frac{\o+x\theta}{\o+x\theta+\bar x} \bigr )$ is
 holomorphic  on  $\CC\setminus (-\infty,0]$. The proof can be concluded 
by  a second application of the same theorems.
\end{proof}

In what follows, we keep the notation $\hat \Pi(\theta, \o)$ for the 
analytic continuation of the Laplace transform.
The next lemma describe the behaviour of $\hat \Pi(\theta,\o)$
near  the origin. Using the  Laplace inversion formula, we  then 
derive from this result the asymptotic behaviour of $\Pi(\theta, t_\o)$
as $t_\o\uto\infty$.

\begin{lemma} \label{cammino}
For any $\theta>0$ let us set
$$
A(\theta) =:\frac{\sin (\pi \a)}{\pi}\int _{\frac{\theta}{\theta+1}} ^1
 u^{-\a} (1-u)^{\a-1} du.
$$
Moreover, let us define
\begin{equation}\label{as}
\cA :=\{r e^{i\phi}\,:\, r\geq 0,\; |\phi|\leq \frac{3}{4}\pi\}.
\end{equation}
Then, for a  suitable positive constant $c>0$
\begin{align}
&  \bigl |\hat\Pi(\theta ,\o)\bigr| \leq c |\o|^{-1},
\qquad \qquad \qquad \forall \o\in\cA:\, |\o|\geq 1 \label{caposud}\\
& \bigl |\hat\Pi(\theta,\o)-A(\theta)/\o\bigr |\leq c |\o|^{-\a},
\quad \;\;\forall \o\in\cA:\, |\o|\leq 1\label{caponord} .
\end{align}
\end{lemma}
\begin{proof}
The first estimate (\ref{caposud}) follows trivially from
 (\ref{davide}) and (\ref{palline}). Let us prove  (\ref{caponord})
 for 
$\o\in \cA$ and $|\o|\leq 1$.

In what follows,  $c_0,c_1,\dots$ denote some suitable
positive constants depending 
only on $\theta$. Moreover,  given $z\in\CC$ we denote by $\int _0^z $ 
and $\int_z^\infty$ the integrals over the paths 
$\{ sz\,:\, 0\leq s\leq 1\}$ and $\{sz\,:\, s\geq 1\}$, respectively. 
We extend the functions  $z^{-\a}$ and $z^{\a-1}$, 
defined on $(0,\infty)$, to $\CC \setminus (-\infty,0]$ by analytic
  continuation.  Then (\ref{davide}) implies
\begin{equation}\label{roma21}
\begin{split}
\o \hat\Pi(\theta, \o) &=
\int_0 ^{1/\o} x^{\a-1}
\Bigl(
[1+ x(1+\theta)]\,[1+x\theta]
\int_0^{1/\o}  \frac{y^{\a-1}}{1+x\theta +y }dy \Bigr)^{-1} dx \\
& = \int_0 ^{1/\o} x^{\a-1}
\Bigl(
[1+ x(1+\theta)]\,[1+x\theta] ^\a
\int_0^{\frac{1}{\o(1+x\theta)}}  \frac{y^{\a-1}}{1+y }dy \Bigr)^{-1} dx.
\end{split}
\end{equation}
Let us define 
$$
\cB :=\{ (\o,x)\in \CC^2 \;\text{ s.t. }
\;\o\in \cA,\;  x=\frac{s}{\o} \text{
for some } s: 0\leq s\leq 1 \}.
$$
Since  
$\bigl( \o(1+x\theta)\bigr)^{-1} \in \cA\cap \{z\,:\, |z|\geq c_0\}
$, we obtain
\begin{align}
& \Bigl|\int_0^{\frac{1}{\o(1+x\theta)}}  \frac{y^{\a-1}}{1+y }dy \Bigr|
\geq c_1, \label{nero1}  \\
&  \Bigl|\int_{\frac{1}{\o(1+x\theta)}}^\infty  \frac{y^{\a-1}}{1+y }dy
 \Bigr|\leq c_2 \,  |\o(1+x\theta) |^{1-\a} .\label{nero2}
\end{align}
Let $B(\theta)$ be defined as 
\begin{equation}\label{roma24}
B(\theta):= \int_0 ^{\infty} \frac{y^{\a-1}}{1+y }dy=
\int _0 ^1 u^{-\a} (1-u)^{\a-1} du=\frac{\pi}{\sin (\pi \a)}.
\end{equation}
(note that the above second identity  follows from
 the change of variable $u=y(1+y)^{-1}$, while the last one is well 
known in the theory of the  Gamma function).
By means of  (\ref{nero1}) and (\ref{nero2}) we obtain
\begin{equation}
\begin{split}
 \Bigl|
\o \hat\Pi(\theta, \o) -\frac{1}{B(\theta)} \int _0 ^{\frac{1}{\o}} &
\frac{x^{\a-1}dx}{ \bigl(1+ x(1+\theta)\bigr)\,\bigl(1+x\theta\bigr)^\a }
\Bigr| \leq \\
&  |\o|^{1-\a} 
\int _0 ^{\frac{1}{\o}} 
\frac{|x|^{\a-1}\, d|x|}{\bigl|1+ x(1+\theta)\bigr|\,\bigl|1+x\theta\bigr|^{2\a-1}
  }\leq c_3   |\o|^{1-\a} .
\end{split}
\end{equation}
Since
$$ 
 \Bigl|  \int _{\frac{1}{\o}} ^\infty
\frac{x^{\a-1}dx}{ \bigl(1+ x(1+\theta)\bigr)\,\bigl(1+x\theta\bigr)^\a }
\Bigr| \leq c_4 |\o| 
$$
and, using analyticity  and integrability of the
singularities around $z=0$ and $z=\infty$,
$$ 
\int _{s\o:\, s\geq 0 }  
\frac{x^{\a-1}dx}{ \bigl(1+ x(1+\theta)\bigr)\,\bigl(1+x\theta\bigr)^\a }=
\int  _0 ^\infty  
\frac{x^{\a-1}dx}{ \bigl(1+ x(1+\theta)\bigr)\,\bigl(1+x\theta\bigr) ^\a },
$$
we get
$$
\Bigl | \o \hat \Pi (\theta, \p) - \frac{1}{B(\theta)}
 \int_0 ^{\infty} \frac{x^{\a-1} }{
(1+ x(1+\theta))(1+x\theta)^\a }dx \Bigr| \leq c_5 |\o|^{1-\a}.
$$
Using the chance of variables  
$v=x^{-1}+\theta$ and $u= v(1+v)^{-1}$, we obtain
$$
 \int_0 ^{\infty} \frac{x^{\a-1} }{
(1+ x(1+\theta))(1+x\theta)^\a }dx=
\int _{\frac{\theta}{\theta+1}} ^1
 u^{-\a} (1-u)^{\a-1} du
$$
which implies the assertion of the lemma.
\end{proof}

Lemma \ref{cammino} and Proposition \ref{LTI} allow us to conclude
the proof of the ageing behaviour of $\Pi_N(t,t_w) $:
\begin{prop}\label{tatoo}
 For almost all energy landscapes $\underline E$,  given $\theta >0$
  \begin{equation}\label{tatoorom}
    \lim_{t_w\uto\infty}\lim_{N\uto\infty} \Pi_N (\theta t_w,t_w)=
  \frac{\sin (\pi \a)}{\pi}  \int_{\frac \theta {1+\theta }}^1 u^{-\a}(1-u)^{\a-1}du.
  \end{equation}
\end{prop}

%
%

\subsection{Visiting deeper and deeper traps}\label{step2}\\

In this section we use the integral representation
(\ref{caporale}) in order to study the probability that the system
at time $t$  is in a deep trap, i.e. in a state  with large
waiting time. In Proposition \ref{bruna} we first prove that  the probability to be in a site
with waiting time smaller than $O(1)$ decays as $t^{\a-1}$, thus
implying the ageing behaviour of other correlation functions
described in subsection \ref{step3}. In the second part, we will
investigate the random variable $t x_N(t)$ and show that, for
almost all $\underline E$,  it has a weak limit as $N\uto\infty$ and
then $t\uto\infty$. As consequence, with high probability at time
$t$ the system is in a state  of waiting
time $O(t)$ as stated in Proposition \ref{nicaragua}.\\

Reasoning as
in the proof of Proposition   (\ref{schneider}), we can prove,
 for almost all energy landscapes $\underline E$, that, given a
 function $h$ on $[0,1]$ that can be uniformly approximated 
by piecewise $C^{1}$ functions, 
\begin{equation}\label{wyler}
H(t):=\lim_{N\uto \infty}\EE_N\bigl(\,h( x_N(t)\,\bigr) =
\frac{1}{2\pi i}\int_\g \frac{e^{- t \l}}{\l} \frac{ \int _0 ^1
\frac{h(x) }{\l-x} x^{\a-1}dx}{\int _0 ^1 \frac{1}{\l-x}x^{\a-1}
dx} d\l \quad \forall t>0,
\end{equation}
where $\g$ is a positive oriented loop around $[0,1]$.

Since $H(t)$ is a bounded function, the Laplace integral $ \hat
H(\o):=\int _0^\infty  H(t) e^{-\o t} dt $ is absolutely
convergent when $\Re (\o)>0$. By the same arguments we used 
to derive (\ref{davide}), it is simple to deduce from the integral
representation (\ref{wyler})  that
\begin{equation}\label{miele}
 \hat H(\o) =\frac{1}{\o}
\frac{ \int _0 ^1 \frac{h(x) }{\o+x} x^{\a-1}dx}{\int _0 ^1 \frac{1}{\o+x}x^{\a-1} dx}
\end{equation}
In the following Proposition  we concentrate on  the case
 $h(x):=\II_{x\geq \d}$. By (\ref{miele}), we can give precise 
information on the  asymptotic behaviour of the probability to be
 at time $t$ in a site with waiting time
 smaller than $1/\d$:

\begin{prop}\label{bruna}
Let
\begin{equation}
B(\d):=\frac{  \int _{\d} ^1 x^{\a-2} dx}{ \int_0^\infty
\frac{x^{\a-1}}{1+x} dx },\qquad c(\a):= 
\int_0 ^\infty y^{\a-1}e^{-y} dy
\end{equation}
Then, for almost all energy landscapes $\underline E$,
\begin{equation}\label{zivago}
\lim_{s\uto\infty} s^{1-\a}\lim _{N\uto \infty} \PP_N 
\bigl (x_N(s)>\d\bigr)=
 B(\d)/c(\a) .
\end{equation}
\end{prop}

Finally, we show that with high probability at time $t$ 
 the system is in a trap of depth of order $O(t)$.  In particular,
 the random variables $\,tx_N(t)$ converge weakly to a nonnegative random variable  as $N\uto\infty$ and then  $t\uto\infty$ a.s.. This result corresponds to the  convergence of expectation of bounded continuous functions and due to Lemma
\ref{pigiama} such a  convergence can be extended to the larger class of bounded
piecewise continuous functions, which is more suitable in order to investigate
the phenomenon of visiting deeper and deeper traps:
\begin{prop}\label{nicaragua}
Let $Z$ be the only random variable with range in $(0,\infty)$ having
Laplace transform
$$ \EE\bigl ( e^{-\theta Z} \bigr )= \frac{\sin (\pi \a)}{\pi}
\int _{\frac{\theta}{\theta+1}} ^1
 u^{-\a} (1-u)^{\a-1} du.
$$
Then, for almost all energy landscape $\underline E$, given a
bounded piecewise continuous  function $h$ on $(0,\infty)$,
\begin{equation}\label{tt667}
\lim_{t\uto\infty} \lim _{N\uto\infty} \EE_N ( h(t x_N(t)) )=
\lim_{t\uto\infty}
\frac{1}{2\pi i}\int_\g
\frac{e^{- t \l}}{\l}
\frac{ \int_0^1 \frac{h(xt)}{\l-x}x^{\a-1}dx  }{\int_0^1 \frac{1}{\l-x}
x^{\a-1}dx } d\l
=
 \EE( h (Z))
\end{equation}
In particular, for almost all energy landscapes $\underline E$,
$$
\lim _{t\uto\infty} \lim_{N\uto\infty}
 \PP ( \frac{\t _N (t)}{t} \geq u ) =
\PP ( Z \leq u^{-1} ),\qquad \forall u >0.
$$
\end{prop}

\vspace{.3cm}

\noindent
\emph{Proof of Proposition \ref{bruna}}.
 We have to prove that
$\lim _{s\uto\infty} s^{1-\a} H(s)= B(\d)/c(\a)$ where
 $H$ is given by (\ref{wyler}) with $h(x):=\II_{x\geq \d}$. 
As in the proof of Lemma \ref{pisa} we can show  that the r.h.s. of
 (\ref{miele}) is well defined and holomorphic on
$\CC\setminus (-\infty, 0] $.  We keep the notation 
$\hat H$ for this extended function.
 By the change of variables $x=\o y$, we get
\begin{equation}
\int_0 ^1 \frac{x^{\a-1}}{\o+x}dx = \o^{\a-1} \int _{\g_w} 
\frac{y^{\a-1} }{1+y} dy
\end{equation}
where $\g_\o$ is the oriented path
$\{s/\o\}_{ 0\leq s\leq 1}$. Let $\hat \g_\o$ be the path
$\{s/\o\}_{ s\geq 0 }$.
 By analyticity and integrability of the singularities
 at $z=0, z=\infty$, we have
$$
\int _{\hat \g_\o}\frac{y^{\a-1} }{1+y} dy   = 
\int  _0 ^\infty \frac{y^{\a-1} }{1+y} dy
$$
Let us define
$\cA:=   \{ re^{i\theta}\,:\, 0<r<\infty,\; |\theta|\leq \frac{3}{4}\pi \}
$. Then, for a suitable constant $c_1$,
$$
\bigl| \int _{\hat \g_\o\setminus \g_\o}\frac{y^{\a-1} }{1+y} dy\bigr|
\leq c_1  |\o |^{1-\a},\qquad
\forall \o\in \cA\,:\,|\o|\leq 1,
$$
implying
\begin{equation}\label{natale1}
\int_0^1\frac{x^{\a-1}}{\o+x}dx =\o^{\a-1}\Bigl(
 \int  _0 ^\infty \frac{y^{\a-1} }{1+y} dy + O(  |\o |^{1-\a} ) \,\Bigr),
\end{equation}
where $A=B+O(1/N)$ is understood to mean that there exists $C<\infty$
 such that  $|A-B|\leq C/N$.
Trivially,
\begin{equation}\label{natale2}
\int_\d^1\frac{x^{\a-1}}{\o+x}dx=\bigl(1+O(|\o|)\bigl)
\int_\d ^1  x^{\a-2}dx.
\end{equation}
Let us note that  the estimate of error terms   in 
(\ref{natale1}) and in (\ref{natale2}) is  uniform in $\o\in \cA$, $|\o|\leq 1$. 
Then, from  (\ref{miele}), (\ref{natale1}),
 (\ref{natale2}) we get
\begin{equation}\label{pasqua1}
|\o^\a \hat H(\o) - B(\d)/c(\a)|\leq c_2|\o |^{1-\a},
\quad  \forall \o\in \cA\,:\, |\o|\leq 1.
\end{equation}
Since trivially $|\hat H(\o)| \leq c_3 |\o|^{-1}$ for 
$\o \in \cA $ with $|\o|>1$, the assertion of the proposition 
follows from Proposition \ref{LTI}.
\qed

\vspace{0.5cm}

\emph{Proof of Proposition \ref{nicaragua}}.
As discussed before (\ref{wyler}), one can show,  for almost all energy
landscapes $\underline E$, that, given a piecewise continuous function $h$
on $(0,\infty)$,
\begin{equation*}
\Phi_t(h):= \lim _{N\uto\infty} \EE_N ( h(t x_N(t)) )=
\frac{1}{2\pi i}\int_\g \frac{e^{- t \l}}{\l} \frac{ \int_0^1
\frac{h(xt)}{\l-x}x^{\a-1}dx  }{\int_0^1 \frac{1}{\l-x} x^{\a-1}dx
} d\l\qquad \forall t\geq 0,
\end{equation*}
where $\g$ is a positive oriented loop around $[0,1]$.
Note that $\Phi _t$ defines a positive linear functional on the
space of continuous functions on $(0,\infty)$, decaying at $\infty$
and  satisfying  $\Phi_t(1)=1$. Therefore,  the  Riesz--Markov 
representation theorem (see
Theorem IV.18 in \cite{RS}) implies that  $\Phi_t(h)= \mu_t (h)$,
 for a unique Borel probability measure $\mu_t$ on $[0,\infty)$. In particular,
there exists a  random variable $Z_t$ on $(0,\infty)$ such that 
$$
\lim _{N\uto\infty} tx_N(t) \rightarrow Z_t \text{ weakly }, \forall t>0
\qquad \text{a.s.} .
$$
 If we take $h(t)=e^{-t \theta}$, then 
$\Phi_t(h)=\mu_t(h)= \Pi(\theta t, t)$, with $\Pi$ defined as in
(\ref{iacopo11}). That means that   $\Pi(\theta t, t)$ is the Laplace transform
of $Z_t$. As proved in subsection \ref{step1},
$$ 
\lim_{t\uto\infty}\Pi(\theta t, t)  = \frac{\sin (\pi \a)}{\pi}
\int _{\frac{\theta}{\theta+1}} ^1
 u^{-\a} (1-u)^{\a-1} du:=f(\theta) .
$$
We state that $f(\theta)$  is the Laplace transform of a random
variable $Z$ with range in $[0,\infty)$. To this aim we apply the
criterion given by Theorem 1, Section XIII.4   in \cite{Feller}.
By (\ref{roma24}), $f(0)=1$. Moreover, $f^{(1)} (\theta) =
- \frac{\sin (\pi \a)}{\pi}\theta^{-\a}(1+\theta)^{-1}$ thus
implying (by trivial computations) that $(-1)^k f^{(k)}(\theta)
\geq 0$. This completes the proof of our statement.

Since the Laplace transform of $Z_t$ converges to the Laplace 
transform of $Z$ as $t$ goes to $\infty$, we have that $Z_t$ 
converges weakly to $Z$, thus implying 
 the limit (\ref{tt667}) whenever $h$
is a bounded continuous functions on $(0,\infty)$.  
Finally, due to  Theorem 5.2 in \cite{billy}, the limit remains valid
if $h$ is a bounded measurable function whose set of discontinuity points has
zero measure w.r.t. the distribution of $Z$. Therefore, 
   Lemma \ref{pigiama} allows to prove  (\ref{tt667}) for $h$ bounded and piecewise continuous.
\qed

\vspace{0.2cm}

\begin{lemma}\label{pigiama}
The distribution function $F(z):=\PP(Z\leq z)$ of the positive
random variable $Z$ is  continuous.
\end{lemma}
\begin{proof}
Trivially  $F$ is increasing and right continuous. Therefore, it has a 
countable set of points
of discontinuity. Moreover, by the Laplace inversion formula
 (see \cite{Feller}, XIII.4), if $x$ is a point of
continuity, then
$$ 
F(x)=\lim _{a\rightarrow \infty} \sum _{n\leq a x} \frac{ (-a)^n}{n!}
f^{(n)} (a).
$$
Given $s=0,1,2,\dots $ and $\g>0$ let  $c_s(\g)>0$ be such that
$D_a^{s} a^{-\g}=  (-1)^s c_s(\g) a^{-\g-s}$. Then the Leibniz formula 
implies
$$
(-1)^n D_a^{n} \bigl( a^{-\a}(1+a)^{-1}) =\sum_{s=0}^n c_s(\a) c_{n-s}(1) 
a^{-\a-s} (1+a)^{-1-n+s} \leq (-1)^n D_a^{n} a^{-\a-1}.
$$
Since $f^{(1)} (a) = -\frac{\sin (\pi \a)}{\pi} 
a^{-\a}(1+a)^{-1}$, the above estimate implies
$$
(-1)^n f^{(n)} (a)=|f^{(n)} (a)|\leq  \prod _{k=1}^{n-1}
(k+\a) \theta ^{-n-\a}, \qquad \forall n\geq 2.
$$
In particular, given two points of continuity $0<x<z$ we have
\begin{equation}\label{gnummi100}
F(z)-F(x) \leq  \limsup  _{a\rightarrow \infty} a^{-\a} \sum _{ax<n\leq a z}
\frac{1}{n} \prod _{k=1}^{n-1} (1+\frac{\a}{k}).
\end{equation}
One can prove  that the sequence
$\prod _{k=1}^{n-1} e^{-\frac{\a}{n}}(1+\frac{\a}{k})$ is convergent (see
\cite{ahl}, Chapter 5, Section 2.4). Let us denote its limit by
 $c_\a$ and let $\g$ be  Euler's constant
$$
\gamma:=\lim _{n\uto\infty}\Bigl(1+\frac{1}{2}+\frac{1}{3} +\cdots
+ \frac{1}{n}-\log n \Bigr).
$$
Then we can write
\begin{equation}\label{abudabi}
 \frac{1}{n} \prod _{k=1}^{n-1} (1+\frac{\a}{k})=
 e^{\a\bigl(
1+\frac{1}{2}+\cdots+\frac{1}{n} -\log (n-1) \bigr )}
\frac{(n-1)^\a}{n}\prod _{k=1}^{n-1} e^{-\frac{\a}{k}}(1+\frac{\a}{k}).
\end{equation}
In particular, in (\ref{abudabi}) we can substitute 
$ \prod _{k=1}^{n-1} e^{-\frac{\a}{k} } (1+\frac{\a}{k})$
with $c_\a$  with an error term in (\ref{gnummi100}) bounded by
$$
\text{const.}\, a^{-\a}(az-ax) (ax)^{-1+\a}\Bigl| \prod _{k=1}^{n-1}  e^{-\frac{\a}{k} } (1+\frac{\a}{k})-c_\a\Bigr|
\leq c(x,z)  \Bigl|\prod _{k=1}^{n-1} e^{-\frac{\a}{k} }  (1+\frac{\a}{k})-c_\a\Bigr|
$$
which is negligible as $a\uto\infty$. Therefore
\begin{equation}\label{dubai}
F(z)-F(x) \leq  \limsup  _{a\rightarrow \infty} c_\a a^{-\a} e^{\g \a}  
\sum _{ax<n\leq a z}
(n-1)^{-1+\a}\leq  c' (z^\a-x^\a)
\end{equation}
for a suitable positive constant $c'$.
Since (\ref{dubai}) is valid almost everywhere and $F$ is monotonic, we have
that $F$ is continuous.
\end{proof}

%
%

\section{Other correlation functions}\label{step3}

In this section we study the asymptotic behaviour of some other
time-time correlation functions $\P^{(1)}_N(t,t_w)$,
$\P^{(2)}_N(t,t_w)$ for which  deep traps play a special role.
This section is mainly a preparation of what is to follow in the second
part of the paper.

Given $\d>0$, we define 
 the set of sites with small waiting time as
 $D_N:=\{i\,:\, x_i\geq \d, \; i=1,\dots , N\} $. 
Moreover, we set
\begin{align}
& \P^{(1)}_N(t,t_w):=  
\PP_N \bigl(\, Y_N (u)\in D_N  \quad \forall u\in (t_w,t_w+t]\text{ s.t. }
 Y_N(u)\not = Y_N(u^-)\;\bigr)\\
& \P^{(2)}_N(t,t_w):= 
  \PP_N \bigl(\, Y_N(u) \in D_N \cup\{ Y_N(t_w) \}
\quad \forall u\in (t_w,t_w+t] \text{ s.t. } x_N(u)\not = x_N(u^-)\;\bigr).
\end{align} 

Given a subset $A\subset \cS_N$, $i\in A$ and $s>0$, let $\varphi _{N,A} (i,s) $ be defined as 
$$ 
 \varphi _{N,A}(i,s):= \PP_N\bigl( Y_N(u)\in A \;\forall u\in [0,s]\,|\, 
Y_N(0)=i\bigr),
$$
then
\begin{align}
& \P^{(1)}_N (t,t_w)= \P_N (t,t_w) +\sum_{j=1}^N \PP_N(Y_N (t_w)=j)
\int_0^t ds\,\frac{x_j e^{-s  x_j }}{N}\sum _{i\in D_N}
\varphi _{N, D_N} (i, t-s)  \label{cuernavaca1},\\ 
& \P^{(2)}_N (t,t_w)= \sum_{j=1}^N \PP_N(x_N(t_w)=x_j)\, \varphi
_{N, D_N\cup \{j\} } (j,t)\label{cuerna1}
\end{align}
where the first identity can be derived by 
 conditioning on the  first jump performed  after the waiting time $t_w$ and
 by  recalling the  following realization of the dynamics:  after arriving at the state $i$, the system waits an exponential time with parameter $x_i$ and after that it jumps to a site in $\cS_N$ with uniform probability.

The following proposition  is mainly a consequence of the  phenomenon
of visiting deep traps with higher and higher probability.
 To this aim recall Proposition \ref{bruna} or simply, as a
 consequence of Proposition \ref{nicaragua}, that
\begin{equation}\label{schiele}
\lim _{t\uto\infty} \lim _{N\uto\infty} \PP_N( x_N(t)>\ep)=0,\qquad
\forall \ep>0 .
\end{equation}

\begin{prop}\label{rahel}
For almost all $\underline{x}$
\begin{equation}\label{coppito102}
\lim _{t_w\uto\infty} \sup_{t\geq 0} 
| \P^{(i)}_N(t,t_w)-  \P_N (t,t_w)\bigr|=0\qquad \text{\rm for } i=1,2
\end{equation}
\end{prop}
\begin{proof}  We consider first the case $i=1$. 

We claim that 
 for any $u>0$ and $i\in D_N$,
\begin{equation}\label{fata}
 \varphi _{N,D_N} (i,u) \leq   
\exp \Bigl(-\d\,u\bigl(1-\frac{|D_N|}{N} \bigr) \Bigr). 
\end{equation}
In order to prove such a bound, we introduce a new random walk $ Y^\ast_N(t)$
having generator $\LL^\ast$ defined as the r.h.s. of (\ref{perleaiporci}), where $x_i$ is replaced by $\d$ if $i\in D_N$. By a simple coupling argument 
one gets
$$  \varphi _{N,D_N}(i,u)\leq \varphi ^\ast _{N,D_N}(i,u)  $$
where the function  $\varphi ^\ast _{N,D_N}$ is the analogous of 
$ \varphi _{N,D_N}(i,u)$ for the random walk $ Y^\ast_N(t)$. At this point, 
it is enough to observe that  $\varphi ^\ast _{N,D_N}$ equals  
the r.h.s. of (\ref{fata}).

Now fix $\ep>0$. Then, due to (\ref{cuernavaca1}) and (\ref{fata}),
\begin{equation}
\begin{split}
 \bigl|
\P^{(1)}_N & (t,t_w) - \P_N (t,t_w)
\bigr| \\
& \leq  \PP_N( x_N (t_w)\geq \ep)+
  \sum_{j\,:\, x_j <\ep } \PP_N(Y_N (t_w)=j)
  \frac{x_j}{N}\sum _{i\in D_N} \int_0^t 
  \varphi _{N,D_N} (i,t-s) \,ds \\
&  \leq \PP_N( x_N (t_w)\geq \ep)+ \ep 
\int _0 ^t  e^{-\d \,u  \bigl(1-\frac{|D_N|}{N} \bigr) }du.
\end{split}
\end{equation}
By the law  of large numbers,  
$$ 
\lim _{N\uto\infty} 
\int _0 ^\infty  e^{-\d\,u \bigl(1-\frac{|D_N|}{N} \bigr) }
<\infty ,\qquad \text{a.s.} . 
$$
The proposition  now  follows  from  the fact that  $\ep$ is arbitrary
and from (\ref{schiele}).

To deal with case (ii), one proceeds in essentially the same way, 
decomposing the path of the process at its returns to the point $x_i$,
and summing over the number of these returns. One finds easily that the
case when the process does not leave $x_i$ for the entire period $t$
dominates, leading to the assertion of the proposition. We leave the
details to the reader.

\end{proof}

%
%

\section{The REM--like trap model on a Poisson point process.}\label{annan}

In this section we consider a slightly different formulation of the
REM like trap model that betrays more directly its connection to the 
 REM dynamics (see \cite{BBG1,BBG2}) and that offers a
somewhat more natural insight in the r\^ole of time scales in the
analysis of ageing systems.
Let us consider a Poisson point process  
$\cP =\sum_i \d_{E_i} $  on $\RR$ with
intensity  measure  $\a e^{-\a E}dE$, where $0<\a<1$. 
Note that such processes arise naturally as the extremal process of
sequences of random variables.
Trivially, 
a.s. the support of $\cP$ is an infinite set of points, 
whose maximum  is almost surely bounded from above. 
Thus  we can   label the points in the support of 
$\cP$ in decreasing order: $E_1>E_2>\dots$. 
Then, the energy landscape $\underline E$ is defined as
$\underline E=(E_1, E_2, \dots )$. We want to define a random process
on the support of this point process that jumps ``uniformly'' from any
point to any other point in the support. To do this, we need to
introduce a cut--off.  Here we  
fix an energy threshold $E$ and set
$$
N_E=\max\{i\,:\, E_i \geq E  \}.
$$
Note that $N_E$ is a Poisson random  
variable with expectation $e^{-\a E}$. Moreover, 
the probability that $N_E\geq 1$ can be made as small as desired 
  when $E$ is chosen  small enough, 
as we assume in what follows.\\

Let  $\cG _E=(\cS_E,\cE_E)$ be the  graph with
$$
\cS_E:=\{ 1, 2, \dots,  N_E\},\qquad
\cE_E:=\{ \,\{i,j\}\,:\, i\not =j\in \cS_E \}.
$$
Since here we want to investigate the effect of \emph{time
rescaling}, we introduce a time unit $\t_0=e^{E_0}$. Then, the
  continuous--time random walk $Y_E(t)$  is the  random walk on
  $\cG_E$ 
having uniform initial distribution and such that, after arriving at 
site $i\in \cE_E$, it  waits an exponential time with mean
$\frac{N_E}{N_E-1}e^{E_i}/\t_0 $ and then jumps with uniform
probability to a different site of $\cS_E$. In particular, the
Markov  generator $\LL_E$ for the above defined random walk is
given by $\LL_N$ in  (\ref{perleaiporci})  with $N:=N_E$ and
$x_i:= \t_0 e^{-E_i}$ (since for $E\ll 0$, $\frac{N_E}{N_E-1}\sim
1$ when referring to waiting time we disregard the coefficient
$\frac{N_E}{N_E-1}$ as in section \ref{inizio}). Although $Y_E(t)$
depends on $\t_0$, our notation does not  refer to such a
dependence.  In what follows we denote by $\PP_E$ the probability
measure on the paths space determined by $Y_E(\cdot)$, and by
$\EE_E$ the corresponding expectation.

Note that the physical waiting time (the absolute  one) for the
system at state $i$  is given by   $T_i:= e^{E_i}$ while in the
above dynamics the waiting time is $\t_i :=T_i/\t_0$, thus in
agreement with the choice  to consider $\t_0$ as our new time
unit. In what follows we consider, when taking  the thermodynamic
limit $E\dto -\infty$,  three different  kinds of time rescaling:
$\t_0$ fixed,  $\t_0:= e^E $ (that is $E_0=E$)  and    $\t_0  \dto
0$ after $E\dto -\infty$.

As in section \ref{inizio} we are interested in the asymptotic behaviour of
 time--time correlation functions. In particular, let us introduce here the
correlation function 
$$
\Pi_E (t,t_w):=
\PP_E\bigl(\,Y_E(s)=Y_E(t_w),\; \forall s\in [t_w, t_w+t]\bigr).
$$
We will prove that when $\t_0$ is fixed the system exhibits fast
relaxation, thus excluding ageing behaviour (see Proposition
\ref{cita11}). 
At the other extreme,   the  scaling $\t_0= e^{E}$ corresponds to the 
implicit choice made in the standard Bouchaud model considered in the
previous sections. In fact with this choice the system can be thought
of as a \emph{grand canonical} version of the original  
REM--like  trap model and all the results of the previous sections
carry over. 
Finally, we consider the third scaling: $\t_0\dto 0$ after $E\dto -\infty$. 
In Proposition \ref{mamma} we show that when performing such limits
the correlation function $\Pi_E(t,t_w)$ converges to $f(\theta)$ where 
$\theta=t/t_w$ and
$f(\theta)$ denotes the r.h.s. of identity (\ref{tatoorom}), that is  the
limiting behaviour of the correlation function  $\Pi_E(t,t_w)$ is
trivial. At this point a simple consideration is fundamental. If we assume
 that  the physical instruments in the laboratory have sensibility up to the 
time unit $\t_0$, then it is natural to disregard jumps into states with 
physical waiting time 
$T_i=e^{E_i}$  much smaller than $\t_0$. Therefore, a time--time
correlation function much more interesting than $\Pi_E(t,t_w)$ is the
following one, where $\d>0$ is fixed
$$
 \P^{(1)}_E(t,t_w)=  
\PP_E \bigl(\, x_E(u) \geq \d  \quad \forall u\in (t_w,t_w+t]:\;
 x_E(u)\not = x_E(u^-)\;\bigr).
$$
where  $ x_E(t):=x_k$ whenever $Y_E(t)= k. $
In Section \ref{pantani} we prove that $\P^{(1)}_E(t,t_w)$  exhibits
ageing behaviour: $ \P^{(1)}_E(\theta t_w ,t_w)$ converges to the above 
$f(\theta)$ after taking the (ordered) limits $E\dto -\infty$,
$\t_0\dto 0$ and $t_w\uto \infty$.

Finally, in this section we discuss 
 the asymptotic spectral behaviour for  the above time  rescalings.  We will show that  ageing  appears 
 whenever  the
 limiting spectral density has a singularity of order $O(x^{\a-1})$ at $0$.

 Let us recall  some properties  of the Ppp $\sum_i  \d_{x_i}$ with  
intensity measure $\a \t_0^{-\a}x^{\a-1} dx$ on $ (0,\infty)$ which will be 
frequently used below. 
Given $M>0$ the truncated Ppp $\sum _{x_i\leq M}\d_{x_i}$ can be
realized as follows. Let  $n_M$  be a Poisson variable with
expectation $\bigl( \frac{M}{\t_0}\bigr)^\a =\int_0^M
\a\t_0^{-\a}x^{\a-1} dx $ and let $X_i$, $i\in\NN$, be i.i.d.
random variables  on $[0,M]$ with probability distribution
  $p(X)dX= \a M^{-\a} X^{\a-1} dX$. Then
\begin{equation}\label{carusotto}
\sum _{x_i\leq M}\d_{x_i} \sim \sum_{i=1}^{n_M} \d_{X_i},
\end{equation}
 in the sense that the   above point processes have
the same distribution. In particular, by taking $M=\t_0 e^{-E}$,
we get
\begin{equation}\label{camiaf}
\sum _{i\leq N_E}\d_{x_i} \sim \sum_{i=1}^{n^\ast_E} \d_{X_i},
\end{equation}
where $n^\ast_E$ is a Poisson variable with expectation $e^{-\a E}$
and $X_i, i\in \NN,$ are i.i.d. random variables on $[0, \t_0
e^{-E}]$ with probability distribution   $p(X)dX= e^{\a E}\a
\t_0^{-\a} X^{\a-1} dX$.

\noindent
{\bf Notation} It is convenient to introduce the random walks  $x_E(t)$, $\t_E(t)$ defined as
 $$ x_E(t):=x_k, \quad \t_E(t):=\t_k\qquad\text{ if }\, Y_E(t)= k. $$

 We denote by $\g_E$ the positive oriented loop  having support
$$
 \text{supp}(\g_E )= \{x\pm i\,:\, x\in  [-1, \t_0 e^{-E}] \}
\cup \{-1+bi\,:\, |b|\leq 1 \} \cup  \{\t_0 e^{-E}+1 +bi\,:\, |b|\leq 1 \}.
$$
Moreover, we call $\g_\infty$ the infinite open path, oriented
from $\infty+i$ to $\infty -i$, having support
$$
 \text{supp}(\g_\infty )=
 \{x\pm i\,:\, x \geq -1  \}
\cup \{-1+bi\,:\, |b|\leq 1 \} 
$$

Finally, for  given $E$, 
 $0=\l_1 ^{\text{{\tiny $(E)$ }}}  <\l_2  ^{\text{{\tiny $(E)$ }}} 
  <\cdots <\l_{N_E} ^{\text{{\tiny $(E)$ }}} $, are 
 the $N_E$ distinct eigenvalues of the infinitesimal generator $ \LL_E$ (see Proposition \ref{structure1}).

%
\subsection{$\t_0 $ fixed}\label{passo1}\\

Let us first observe  that $ \sum _{i=1}^\infty  \t_i <\infty  $
 for almost all $\underline E$. In fact, since the Ppp  $\sum_i \d_{\t_i}$
has intensity measure $\a\t_0^{\a} \t^{-(1+\a)} d\t$ on $(0,\infty)$, 
\begin{equation*}
\EE(|\{i:\, \t_i \geq 1 \}|) =\int _1 ^\infty \a \t_0^\a
\t^{-(1+\a)} d\t <\infty, \qquad
\EE( \sum _{i:\,\t_i< 1}\t_i)= \int _0 ^1 \a \t_0^\a
\t^{-\a} d\t <\infty.
\end{equation*}

Whenever $\sum _{i=1}^\infty \t_i <\infty$, it is simple to derive 
the asymptotic spectral behaviour of the system from
 Proposition \ref{structure1} and to  show  its fast relaxation, thus 
implying the absence of  ageing:

\begin{prop}
 For almost all $\underline E$,
\begin{equation}
 \lim _{E\dto -\infty} \sum _{j=1} ^{N_E} \d_{\l_j ^{\text{{\tiny $(E)$ }}} 
 } = 
\sum _{j=1}^\infty \d_{\l_j} \qquad \text{ vaguely in } \cM([0,\infty) ) ,
 \end{equation}
where $ \cM ([0,\infty)) $ denotes 
 the space of locally bounded measure on $[0,\infty)$ and
$$ 
\{0=\l_1<\l_2<\l_3<\dots \}
=
\{\l\in \CC\,:\, \sum_{k=1}^\infty \frac{\l}{\l-x_k}=0 \}
$$
\end{prop} 
\begin{proof}
In what follows we assume that  $\sum_{i=1} ^\infty  \t_i<\infty$, which is true a.s.. Then the  function  
$\phi_\infty (\l):= \sum _{k=1}^\infty \frac{\l}{x_k-\l}$ is well defined on
$\CC\setminus \{x_i\,:\, i\geq 1\}$ and
 has non negative zeros $0=\l_1<\l_2<\dots $, such that
$x_{i-1}<\l_i<x_i$
 for any $i>1$. At this point it is enough to 
 show that
$$
\lim _{E\dto -\infty}\l_i^{\text{{\tiny $(E)$ }}}= \l_i , \qquad \forall i=1,2,\dots .
$$
The assertion is trivial for $i=1$. 
Suppose that  $N_E\geq i>1$  and set $\psi_E(\l):=\sum _{k=1}^{N_E} 
\frac{1}{x_k-\l}$.
 Due to  Proposition \ref{structure1}, $\l_i^{\text{{\tiny $(E)$ }}}$ is
 the unique zero of $\psi_E(\l)$ in the interval $(x_{i-1},x_i)$. 
In particular,
$$ \psi _E(\l_i)= \psi _E(\l_i)- \psi _E(\l_i ^{\text{{\tiny $(E)$ }}} ) =
\int_{\l_i^{\text{{\tiny $(E)$ }}}}^{\l_i}  \dot{\psi}_E(\l)d\l.
$$
Since  $ \dot \psi_E(\l)\geq \frac{1}{(x_i-x_{i-1})^2} $ for all  $\l \in (x_{i-1},x_i)$, we get
$$ 
\bigl|\l_i^{\text{{\tiny $(E)$ }}}-\l_i\bigr |\leq (x_i-x_{i-1})^2 |\psi _E(\l_i)|
$$
and therefore the assertion follows by observing that the identity $\phi_\infty(\l_i)=0$ implies 
$$ 
\bigl|
\psi_E(\l_i)
\bigr|
\leq  \sum_{k=N_E+1} ^\infty 
\frac{1}{x_k- x_i} \downarrow 0 \qquad \text{ as } E\dto -\infty . 
$$
\end{proof}

\begin{prop}\label{cita11}
 For almost all $\underline E$, 
\begin{equation}\label{lapeppa!}
\lim _{t\uto\infty} \lim _{E\dto-\infty}
\PP_E (x_E(t)=x_j)= \frac{ \t_j}{\sum_{k=1}^\infty  \t_k}\quad
\forall j=1,2,\dots
\end{equation}
thus implying
\begin{align}
& \lim_{t_w\uto\infty} \lim_{E\dto -\infty} \Pi_E (\theta t_w,t_w)=0\;\; \forall \theta>0,\label{annibale77}\\
& \lim_{t_w\uto\infty} \lim_{E\dto -\infty} \Pi_E ( t,t_w)=
\frac{ \sum _{i=1}^\infty \t_i e^{- x_i t }}{ \sum _{i=1}^\infty \t_i}\;\;
\forall t>0. \label{bruna77}
\end{align}
\end{prop}
\begin{proof}
In what follows we assume that $\underline E$ satisfies
 $\sum _i \t_i <\infty$. 
By setting $h(x)= \II_{ x=x_j}$ in (\ref{caporale}) we get
 the integral representation
\begin{equation}
\PP_E (x_E (t)=x_j)= \frac{1}{2\pi i} \int _{\g_E} \frac{e^{-\l t}}{\l(x_j-\l)}
\bigl( \sum_{k=1}^{N_E}\frac{1}{x_k-\l} \bigr)^{-1}d\l.
\end{equation}
By applying  the residue theorem (see the  arguments used in order to
derive (\ref{davide})),
it is simple to compute the Laplace transform $\hat F_E(\o)=\int _0 ^\infty
\PP_E (x_E (t)=x_j) e^{-\o t } dt$ for $\Re(\o)>0$:
\begin{equation}
\hat F_E(\o)= \Bigl( \o (\o+x_j) \sum_{k=1}^{N_E} \frac{1}{\o+x_k} \Bigr)^{-1}
\end{equation} 
We note that, almost surely, there exists $c>0$ such that
\begin{equation}\label{compleannojiri}
  \bigl|\hat F_E (\o)\bigr|\leq c\, \frac{1}{|\o|},\qquad \forall E, \quad 
\forall \o\in \cA :=
 \{ re^{i\theta}\,:\, 0<r<\infty,\; |\theta|\leq \frac{3}{4}\pi \}
\end{equation}
This follows easily from the bound below where   $\o=a+ib$ and $N$ is  
 any positive integer:
$$
\bigl |
\sum_{j=1}^N \frac{1}{\o+x_j }
\bigr| \geq
\begin{cases}
\frac{1}{|\o+x_1|}  & \text{ if } a\geq 0\\
\frac{b}{(a+x_1)^2 +b^2 } & \text{ if } a<0 .
\end{cases}
$$
Let us now introduce the path $\tilde \gamma$ consisting of the parabolic arcs 
$\{-t\pm i t^2 \,:\, t\geq 1\}$ and the circular arc of radius
$\sqrt{2}$ around the origin connecting (in anti--clockwise way)
$-1-i$ to $-1+i$. The orientation of $\tilde{\gamma}$ is such that
$-1+i$ comes before $-1+i$. Then,
by means of  (\ref{compleannojiri}), the inverse formula of Laplace transform 
and  the dominated convergence theorem, we get
$$
\lim_{E\dto-\infty}
 \PP_E (x_E (t)=x_j)=\int_{\tilde\gamma} e^{t\o}\hat F(\o)d\o,\qquad
F(\o) := \Bigl( \o (\o+x_j) \sum_{k=1}^{\infty} \frac{1}{\o+x_k} \Bigr)^{-1}.
$$
Note that $F(\o)$ is the limit of $F_E(\o)$ as $E\dto \infty$, in
particular it satisfies (\ref{compleannojiri}). Moreover, $F(\o)$ is
the Laplace transform of $ \lim_{E\dto-\infty}  \PP_E (x_E (t)=x_j)$
and trivially
$$
 \bigl| \o F(\o) - \frac{\t_j}{\sum_{k=1}^\infty\t_k} \bigr| \leq c|\o|, \qquad
 \forall |\o|\leq 1:\; \o\in \cA.
$$
At this point (\ref{lapeppa!}) follows from Proposition \ref{LTI}. Moreover,
from (\ref{lapeppa!}) and the identity
$$
 \Pi _E( t, t_w) =\sum _{j=1}^{N_E} \PP_E (\,x_E (t_w)=x_j\,)e^{ -
\frac{N_E-1}{N_E} x_j t }
$$
it is simple to infer (\ref{annibale77}) and (\ref{bruna77}).
\end{proof}

%
%
%

\subsection{$\t_0=e^E $}\label{passo2}\\

Note that by 
 choosing  $\t_0=e^E $, the random variables  $X_1, X_2,\dots $
 introduced in  (\ref{camiaf})  are i.i.d. with distribution  
$p(X)dX= \a   X^{\a-1} dX$ on 
$[0,1]$. Therefore,  due to (\ref{camiaf}), we can think of  $Y_E(t)$ as the
\emph{grand canonical} version of the Bouchaud's REM--like trap
 model. In particular, it exhibits the same asymptotic spectral density
 and the same   ageing behaviour:

\begin{prop}
For almost all $\underline E$, 
$$
\lim _{E\dto -\infty} \frac{1}{N_E} \sum _{j=1}^{N_E}
\d_{\l_j ^{\text{{\tiny $(E)$ }}} 
 } = \a x^{\a-1}dx \qquad \text{ weakly in } \cM([0,1]).
$$
\end{prop}
\begin{proof}
By  approximating continuous functions on $[0,1]$ with 
step functions having  rational values and jumps at rational points, it is enough to prove that, given $0\leq a<b\leq 1$, 
$$
 \lim _{E\dto -\infty}\frac{1}{N_E}
\bigl|\,\{j\,: 1\leq j\leq N_E,\; \l_j^{\text{{\tiny $(E)$ }}}\in [a,b] \}| =
 b^\a-a^\a\qquad \text{a.s.}
$$
We set 
$$
A_E:=\bigl|\{j:\,  1\leq j\leq N_E,\; x_j\in [a,b]\}\bigr| =
\bigl|\{j\,:\, j\geq 1, e^{-E_j}\in [ e^{-E}a, e^{-E}b ]\}\bigr|.
$$
Then, due to Proposition 
\ref{structure1}, we only need to prove that
$$
 \lim _{E\dto -\infty}\frac{A_E}{N_E}=  b^\a-a^\a\qquad \text{a.s.}
$$
To this aim observe that $N_E$ and $A_E$ are Poisson variables with expectation
(or equivalently variance) respectively equals to 
 $e^{-\a E}$ and $e^{-\a E} (b^\a-a^\a)$. Given a positive integer $n$ we set
 $E(n)=-\frac{2}{\a}\ln n $, i.e.
$e^{-\a E(n)} =n ^2$.
It is simple to derive from Chebyshev inequality and 
Borel--Cantelli lemma that
$$
\lim _{n\uto\infty}
 \frac{N_{E(n)}}{e^{-\a E(n)}}=1,\qquad
\lim _{n\uto\infty}
 \frac{A_{E(n)} }{e^{-\a E(n)}}=b^\a-a^\a
\qquad 
\text{a.s.}
$$
By a simple argument based on monotonicity, one can extend the first limit to 
$\lim _{E\dto-\infty} \frac{N_E}{e^{-\a E}}=1$ a.s. In order to extend the second limit to general $E$ we observe  that, whenever  $E(n+1)<E\leq E(n)$,
$$
\bigl | A_E-A_{E(n)}\bigr|\leq \bigl| \bigl\{\,j\,:\, e^{-E_j}
\in [a \,e^{-E(n)}, a\,e ^{-E(n+1)}]\cup
  [b \,e^{-E(n)}, b\,e ^{-E(n+1)}]
\,\bigr\}\bigr|.
$$
Since the r.h.s. is a Poisson variable with expectation of order $O(n)$, by means of Chebyshev inequality and Borel--Cantelli lemma we obtain 
$$ 
\lim _{n\uto \infty }
\sup _{E(n+1)<E\leq E(n)} \frac{ | A_E-A_{E(n)}\bigr|}{ e^{-\a E(n)} }=0\qquad \text{a.s.}
$$ 
thus allowing to prove that $\lim _{E\dto-\infty}
 \frac{A_{E} }{e^{-\a E}}= b^\a-a^\a$ a.s. 
\end{proof}

\begin{prop}\label{mihi}
For almost all $\underline E$,
\begin{equation}\label{77volte7}
\lim _{E\dto -\infty} \Pi_E(t,t_w)=
\frac{1}{2\pi i}\int_\g \frac{e^{- t_w \l}}{\l}
\frac{\int_0 ^1 \frac{e^{-xt}}{\l-x}x^{\a-1}dx }{
\int_0 ^1 \frac{1}{\l-x} x^{\a-1}dx } d\l,\qquad \forall t,t_w
\end{equation}
In particular, for almost  all $\underline E$, given $\theta>0$
 \begin{equation}
    \lim_{t_w\uto\infty}\lim_{E\dto - \infty} \Pi_E (\theta t_w,t_w)=
  \frac{\sin (\pi \a)}{\pi}  \int_{\frac \theta {1+\theta }}^1 u^{-\a}(1-
u)^{\a-1}du.
 \end{equation}
\end{prop}

\begin{proof}
Our starting point is (\ref{camiaf}) and  the following inequality, valid for
 any  bounded function $f$ with
 $\EE(f(X_i))=0$:
$$
\PP\bigl(\,\bigl |\Av_{j=1}^k f(X_j)\bigr |\geq \d\,\bigr)\leq 2
\exp \bigl( -\frac{k\d^2}{ 4 \|f\|_\infty} \bigr), \qquad \forall \d>0, k=1,2,\dots.
$$
In particular,  by conditioning on $n_E^\ast $ (see (\ref{camiaf}) we get
\begin{equation}\label{mimmo}
\PP\bigl(\,\bigl |\Av_{j=1}^{n_E^\ast } f(X_j)\bigr |\geq \d\,\bigr)\leq
2 \exp \Bigl \{
-e^{-\a E}
(1-e^{
-\frac{\d^2}{ 4 \|f\|_\infty}
}  ) 
\Bigr\}
\end{equation}
It is simple to derive (\ref{77volte7}) from 
the above estimate, Borel--Cantelli lemma and  the integral representation
\begin{equation}\label{memoria}
\Pi_E(t,t_w) =
 \frac{1}{2\pi i}\int_{\g} \frac{e^{- t_w \l }}{\l}  \frac{
\Av_{j=1}^{N_E}   \frac{e^{- x_j t } }{ x_j -\l } }{
 \Av_{j=1} ^{N_E} \frac{1}{x_j-\l}} d\l
\end{equation}
where $\g$ is a positive oriented closed path around $[0,1]$ 
(see \ref{integrale}).
Note that the r.h.s. of (\ref{77volte7}) corresponds to the function 
$\Pi(t,t_w)$ introduced in Proposition \ref{schneider}. Therefore, the 
conclusion of the proof follows from Propositions \ref{schneider} and 
\ref{tatoo}.
\end{proof}

%
%
\subsection{$\t_0 \dto 0 $ after $E\dto -\infty $}\label{passo3}\\

 In this scaling regime, we show that  the vague limit of the suitably rescaled spectral density is given by  the measure $\a x^{\a-1}dx$ on $[0,\infty)$ and 
 we recover the ageing property of the
 correlation function as before. What is more, however, is that due to
 the fact that we are effectively already at 'infinite times' on the
 microscopic scale, we get a pure ageing function even before
 taking $t$ and $t_w$ to infinity:
\begin{prop}\label{sforzetto}
For almost all $\underline E$,
$$
\lim _{\t_0\dto 0} \lim _{E\dto -\infty} \t_0^\a  \sum _{j=1}^{N_E}
\d_{\l_j ^{\text{{\tiny $(E)$ }}} 
 } = \a x^{\a-1}dx \qquad \text{ vaguely in } \cM([0,\infty)).
$$
\end{prop}
\begin{prop}\label{mamma}
For almost all energy landscape $\underline{E}$, given positive
$t, t_w$,
\begin{equation}\label{fabri1}
 \lim_{E\dto -\infty }\Pi_E(t,t_w) =
 \frac{1}{2\pi i}\int_{\g_\infty } \frac{e^{- t_w \l }}{\l} \, \frac{
\t_0^\a \sum_{j=1}^{\infty}   \frac{e^{- x_j t } }{ x_j -\l } }{
 \t_0^\a \sum_{j=1 } ^{\infty } \frac{1}{x_j-\l}}\,\, d\l
\end{equation}
and
\begin{equation}\label{noiosetto}
    \lim_{\t_0\dto 0}
\lim_{E\dto - \infty} \Pi_E (t,t_w)=
  \frac{\sin (\pi \a)}{\pi}  \int_{\frac \theta {1+\theta }}^1 u^{-\a}(1-u)^{\a-1}du \qquad \text{ where }\theta=\frac{t}{t_w}.
\end{equation}
\end{prop} 
\noindent
{\bf Remark} The integral in (\ref{fabri1}) exists
due to Lemma \ref{annibale}.\\

Due to Proposition \ref{structure1}, Proposition \ref{sforzetto} follows if one is able to prove that $ \t_0^\a  \sum _{j=1}^{N_E}
\d_{x_j} $ converges vaguely to $\a x^{\a-1}dx $ on $[0,\infty)$ when taking the (ordered) limits $E\dto -\infty$, $\t_0\dto 0$. This is the content of  Lemma
\ref{averppp} below  concerning   the self--average property of Poisson point
processes with finite intensity measure
 (compare it with Lemma 4.16  in \cite{BBG2}). Finally,
the proof of Proposition  \ref{mamma} is based on (and given after)
the technical  Lemmata \ref{averppp},  \ref{annibale}.

\begin{lemma}\label{averppp}
Let $M>0$ and  let $f$ be a bounded continuous function on $[0,M]$.
Then there exists $\d>0$ such that for almost all energy landscape $\underline E$
\begin{equation}\label{giorno}
\bigl|
\t_0^\a \sum _{x_i\leq M } f(x_i) -  \int _0 ^M f(x) \a x^{\a-1} dx 
\bigr| \leq c\, \t_0^\d,\qquad \forall \t_0>0 
\end{equation}
for a suitable positive constant $c$.
\end{lemma}
\begin{proof}
Let $X_1, X_2, \dots $  and $n_M$ be as in (\ref{carusotto}). Due
 to (\ref{carusotto}) and  since $Var(n_M)= \EE(n_M)=  (M/\t_0)^\a$,
$$
\PP\Bigl(\Bigl| \frac{ | \{ j\,:\, x_j\leq M\}|}{ (M/\t_0)^\a } -1 
\Bigr|\geq \ep \Bigr) \leq \bigl(\t_0/M)^\a \ep^{-2}
$$
In particular, given  $\g,s>0$ with $2s-\g\a<-1$,  by means of Borel--Cantelli lemma we obtain that for almost all energy landscape $\underline{E}$ there exists  $c>0$ such that
$$
\Bigl| \frac{ | \{ j\,:\, x_j\leq M\}|}{ (M/\t_0)^\a } -1 
\Bigr|\leq  c\, k^{-s}, \qquad \forall k=1,2,\dots  \text{ where } \t_0:=k^{-\g}.
$$ 
Due to the above estimate,
\begin{equation}\label{trattore1}
\Bigl|  
\t_0^\a \sum _{x_i\leq M} f(x_i) -M^\a {\Av}_{x_i\leq M} f(x_i) 
\Bigr| 
 \leq c k^{-s} \|f\|_\infty,\qquad   \forall k\in 1,2,\dots \text{ where } \t_0:=k^{-\g}
\end{equation}
where ${\Av}_{x_i\leq M}$ denotes the average over the set $\{x_i\leq M\}$. 
As done for (\ref{mimmo}), if $0< \rho<1$,
\begin{multline*}
\PP\Bigl( \Bigl| {\Av}_{x_i\leq M} f(x_i) -
M^{-\a}\int _0 ^M f(x) \a x^{\a-1} dx \Bigr| \geq \rho \Bigr)\\ \leq
2 \exp\bigl\{ -\bigl(\frac{M}{\t_0}\bigr)^\a
 (1 -e^{-c\rho^2 }) \bigr\}\leq 2 e^{- c' \rho^2 \t_0^{-\a}}
\end{multline*}
In particular, by Borel--Cantelli lemma, for almost all $\underline E$,
\begin{equation}\label{trattore2}
\Bigl| {\Av}_{x_i\leq M} f(x_i) -
M^{-\a}\int _0 ^M f(x) \a x^{\a-1} dx \Bigr| \leq c\, k ^{-s} 
,\qquad  \forall k=1,2,\dots \text{ where } \t_0:=k^{-\g}
\end{equation}
if $s$ is chosen small enough. 
At this point (\ref{trattore1}) and (\ref{trattore2}) imply the
assertion of the lemma, if 
 $\t_0 = k^{-\g}$ for some $k=1,2,\dots$ 
The general case $\t_0>0$ follows easily from the uniform continuity of $f$.
\end{proof}

\begin{lemma}\label{annibale}
For almost all energy landscapes $\underline{E}$, there are
positive constants $\t_0^\ast,c_1,c_2$ having the following properties.
If  $\t_0\leq \t_0^\ast$, $N\geq |\{j\,:\,x_j\leq 1\} | $ and 
$\l\in \g_\infty$ (or $\l =a+ib$ with $|b|\leq 1$ and $a\geq x_j+1 $ for all 
$j\leq N$) then
\begin{equation}\label{pisc1}
 \Bigl|\t_0^\a \sum_{j=1 } ^{N}\frac{1}{x_j-\l}\Bigr|\geq
c_1 |\l|^{-2}.
\end{equation}
Moreover, if $\t_0\leq \t_0^\ast$ and  $M\geq 1$, then
\begin{align}
& \t_0^\a \sum _{j=1}^\infty \frac{1}{|x_j-\l|}\leq c_2 |\l|^{\a-1} \ln (1+|\l|)
,\qquad
\text{ if } \l\in \g_\infty  \label{pisc2}\\ 
& \t_0^\a \sum _{x_j\leq M }^\infty \frac{1}{|x_j-\l|}\leq c_2
M^{\a-1}\ln M ,\quad \text{ if } \l\in \gamma_\infty \text{ or }
\Re(\l) = M+1 \label{pisc3} 
\end{align} 
where $ \sum_{x_j\leq M }$ means $\sum_{j\geq 1 : x_j\leq M }$. 
\end{lemma}

\begin{proof}
It is convenient  to introduce the non rescaled Ppp $\sum _i
\d_{y_i} $,  with $y_i:=e^{-E_i}$, having  intensity
measure $\a y^{\a-1} dy$ on $(0,\infty)$. Moreover, we fix here $\b>2$ and
$0<\g<\b/2-1$ and we define
$$
N_n :=|\, \{j\,:\, n^{\b/\a}\leq y_j< (n+1)^{\b/\a} \}\,|.
$$
for $n$ positive integer.
Then a simple application of Borel--Cantelli lemma implies that,
for almost all energy landscapes $\underline{E}$,
\begin{equation}\label{madrileno}
\Bigl| \frac{N_n}{(n+1)^\b-n^\b} -1\Bigr|\leq c\, n^{-\g},\qquad \forall n=1,2,\dots
\end{equation}
 In fact, by Chebyshev inequality and since $N_n$ is a Poisson
variable with $Var(N_n)=\EE(N_n)= (n+1)^\b-n^\b$,
$$
\PP \bigl( | N_n/\EE(N_n)-1\bigr| >n^{-\g} \bigr) \leq c\,
n^{2\g+1 -\b}
$$

Moreover, by Borel--Cantelli lemma and a simple argument based on
monotonicity, it is simple to prove that there exists $\d>0$ such
that, for almost all $\underline E$,
\begin{equation}\label{etica}
\Bigl|\frac{|\{j\,:\, y_j\leq u\} |}{u^\a} -1\Bigr| \leq \k\,
u^{-\d} \qquad \forall u\geq 1.
\end{equation}
for a suitable positive constant $\k$.

In what follows we write  $\l=a+i b$. Then
\begin{equation}\label{fott}
\Bigl|\sum_{j=1 } ^N \frac{1}{x_j-\l}\Bigr| ^2 = \Bigl(\sum _{j=1}
^N \frac{ x_j-a }{ (x_j-a)^2
 +b^2}\Bigr)^2 +
 \Bigl(   \sum _{j=1} ^N \frac{b}{ (x_j-a)^2
 +b^2}\Bigr)^2
\end{equation}

In order to prove (\ref{pisc1}) we assume (\ref{madrileno}) and 
(\ref{etica}) to be valid and  let $0< \t_0\leq \t_0^*\leq 1$ where 
$\t_0^\ast$ is  such
that $u^\a -\k\, u^{-\g}>0$ for all $u\geq 1/\t_0^*$. In particular,
 $\{x_j\,:\, x_j\leq 1\}\not =\emptyset $. By
(\ref{fott}), if  $N\geq |\{j\,:\,x_j\leq 1\} |$ and $\l\in\g_\infty$
$$
\t_0^\a \Bigl|\sum_{j=1 }^N \frac{1}{x_j-\l}\Bigr| \geq
\begin{cases}
\t_0^\a \sum_{x_j\leq 1} \frac{|b|}{(x_j-a)^2+b^2 }\geq c\, |\l|^{-2} \t_0^\a |\{x_j\,:\, x_j\leq 1 \}| &
\text{ if } |b|\geq \frac{1}{2} \\
\t_0^\a \sum_{x_j\leq 1} \frac{x_j+1 }{(x_j+1)^2+b^2 }\geq c
\,\t_0^\a |\{x_j\,:\, x_j\leq 1 \}| & \text{ if }  |b|< \frac{1}{2}
\end{cases}
$$

At this point, (\ref{pisc1}) for $\l\in\g_\infty $ follows from (\ref{etica}). 
The case   $\l=a+i b $, with $|b|\leq 1$ and $a\geq x_j+1 $ for all 
$j\leq N$, can be treated similarly.

It is simple to derive 
 (\ref{pisc2}) and (\ref{pisc3}) from  estimates (\ref{menuchim1}), $\dots$,
(\ref{mendel5}) below valid for almost all $\underline E$:\\

 if $a\leq 100$, $1\leq M$
 and $\l=a+ib \in\g_\infty$ then
\begin{align}
& \t_0^\a \sum _{j=1}^\infty \frac{1}{|x_j-\l|}\leq c, \label{menuchim1}
 \\
& \t_0^\a \sum _{x_j\geq M } \frac{1}{|x_j-\l|}\leq c\,M^{\a-1}; \label{menuchim2} 
\end{align} 

if $\t_0$ is small enough and  $a \geq 100 $ then
\begin{align}
& \t_0^\a \sum _{x_j\leq \frac{a}{2} } \frac{1}{|x_j-a|}\leq c\, a^{\a-1} 
\label{mendel1}\\
& \t_0^\a \sum _{\frac{a}{2}\leq x_j\leq a-1 }
 \frac{1}{|x_j-a|}\leq c\, a^{\a-1} \ln a   \label{mendel2}\\
& \t_0^\a \sum _{a-1\leq  \l \leq a+1 } \frac{1}{|x_j-\l|}\leq c\, a^{\a-1}
  \quad \text{ if } \l=\a+ib \in \g_\infty  \label{mendel3} \\
& \t_0^\a \sum _{a+1\leq x_j\leq 2a  } \frac{1}{|x_j-a|}\leq c\, a^{\a-1} 
\ln a,\label{mendel4} \\
& \t_0^\a\sum _{x_j\geq M }\frac{1}{|x_j-a|}\leq c\,a^{-2\a} M^{\a-1} ,\quad 
\text{ if } M\geq 2a \label{mendel5}
\end{align}

Let   $\l\in \g_\infty$ with $a\leq 100$. Then, due to (\ref{etica}),
$$ 
\t_0^\a \sum _{x_j\leq 1} \frac{1}{|x_j-\l|} \leq c\,\t_0^\a |\{j\,:\, y_j\leq 
\frac{1}{\t_0} \}| \leq c ' 
$$
while, due to (\ref{madrileno}),
\begin{equation}\label{roth1} 
\t_0^\a \sum _{x_j\geq 1} \frac{1}{|x_j-\l|} \leq c\,\t_0^\a \sum _{x_j\geq 1} 
\frac{1}{x_j}  \leq c' \t_0^{\a-1} \sum _{n\geq \lfloor \t_0^{-\a/\b }
\rfloor } n^{\b-1-\frac{\b}{\a}}\leq c''
\end{equation}
thus proving (\ref{menuchim1}). The proof of (\ref{menuchim2}) 
follows the same arguments of (\ref{roth1}).
(\ref{mendel1}) is a simple consequence of (\ref{etica}). The l.h.s. of 
(\ref{mendel3}) can be bounded by $\t_0 ^\a |\{j\,:\, a-1\leq \t_0 y_j \leq
a+1 \}|$ and (\ref{madrileno}) allows to conclude the proof of (\ref{mendel3}).

The proof of (\ref{mendel2}), (\ref{mendel4}) and (\ref{mendel5}) can be  easily  derived from the following estimate. Let $1\leq A\leq B$ with $B\leq a-1$ or $A\geq a+1$, then
 (\ref{madrileno}) implies
\begin{equation*}
\begin{split}
\t_0^\a \sum _{A\leq x_i \leq B} \frac{1}{|x_i-a|} &\leq c\,
\t_0^\a \sum _{n =n_-}^{n=n_+} \frac{ n^{\b-1} }{|\t_0 n^{\b/\a} -1|} \leq c'
\t_0 ^\a \int _u ^v \frac{x^{\b-1} }{|a-\t_0 x^{\b/\a}| } dx \\
& =c' \,a^{-1-\a} 
\int _{u (\frac{\t_0}{a})^{\a/\b} } ^{v (\frac{\t_0}{a})^{\a/\b} }
\frac{y^{\b-1}}{|1-y^{\b/\a}| } dy
\end{split} 
\end{equation*}
where
$n_- =\lfloor (A/\t_0 )^{\a/\b} \rfloor-1 $, $n_+
 =\lfloor (B/\t_0 )^{\a/\b}\rfloor +1$, $u =n_- - 1$, $ v=n_+ +1 $
(we assume $\t_0$ small enough in order to exclude the singular point in the above intervals of  sum and integral).
\end{proof}

\noindent 
{\bf Proof of Proposition \ref{mamma}}. In order to avoid confusion we underline  here the dependence on $\t_0=e^{E_0}$ by writing $\Pi_{E,E_0}(t,t_w)$ in place of $\Pi_E(t,t_w)$.
Our  starting point is given by the integral representation (\ref{integrale}):
\begin{equation}\label{memoria1}
\Pi_{E,E_0} (t,t_w) =
 \frac{1}{2\pi i}\int_{\g_E} \frac{e^{- t_w \l }}{\l}  \frac{
\t_0^\a \sum _{j=1}^{N_E}   \frac{e^{- x_j t } }{ x_j -\l } }{
 \t_0^\a \sum _{j=1} ^{N_E} \frac{1}{x_j-\l}} d\l.
\end{equation}
Let us choose $\underline E$ satisfying Lemma \ref{annibale}.
Then, due to the exponential decaying factor $e^{- t_w \l }$ and to Lemma
\ref{annibale},
 if $\t_0\leq \t_0^\ast$ and $E$ is small enough such that
$ \t_0 e^{-E} \geq 1$, the path integral $\g_E$ in (\ref{memoria1}) can be 
substituted with $\g_\infty$. At this point, (\ref{fabri1}) follows from
 Lemma \ref{annibale} and   the Dominated Convergence Theorem.

In order to prove (\ref{noiosetto}), given a positive integer $M$, we set
$$
g_{M,E_0}(t, t_w):= 
 \frac{1}{2\pi i}\int_{\G_M } \frac{e^{- t_w \l }}{\l} \, \frac{
\t_0^\a \sum_{x_j\leq M }   \frac{e^{- x_j t } }{ x_j -\l } }{
 \t_0^\a \sum_{x_j\leq M  }\frac{1}{x_j-\l}}\,\, d\l
$$
where  $\G_M$ is the positive oriented path having support 
$$ 
\text{supp}(\G_M)=\{\l\in \CC\,:\, |\l-x|=1 \text{ for some } x\in [0,M]\}.
$$
Then, by applying  Lemma \ref{annibale}, whenever  $\t_0\leq
\t_0^\ast$ 
$$
\bigl| \lim _{E\dto-\infty  } \Pi_{E,E_0}(t,t_w) - g_{M,E_0}
(t , t_w) \bigr| \leq c M^{\a-1} \ln M,\qquad \forall M \in \NN_+ .
$$

Let us assume that $\underline E$ satisfies (\ref{giorno}) 
 for all $M\in \NN_+  $ and for $f(x)=\frac{e^{-xt}}{x-\l}$ or
 $f(x)=\frac{1}{x-\l}$, for all $\l$ in a countable dense set of $\G_M$ and for all rational positive $t$. Then, by a chaining argument, we get
$$
\lim_{E_0\dto -\infty} g_{M,E_0} (t,t_w) = g_M(t,t_w),\qquad \forall t,t_w>0
$$
where 
\begin{equation}\label{marcopolo}
g_M(t, t_w) = \frac{1}{2\pi i}\int_{\G_M } \frac{e^{- t_w \l }}{\l} 
\frac{ \int _0 ^M \frac{e^{-x t }}{\l-x} x^{\a-1} dx }{
  \int _0 ^M \frac{1}{\l-x} x^{\a-1} dx } d\l,
\end{equation}
thus implying
$$
 \limsup _{E_0\dto -\infty} \bigl|\,
\lim_{E\dto-\infty} \Pi_{E,E_0} (t,t_w) -g_M (t,t_w) \bigr| \leq
 M^{\a-1} \ln M,\qquad \forall M \in \NN_+ .
$$

At this point let us observe (see the proof of Lemma \ref{annibale}) that there exist $c,c'>0$ such that for all $M\in\NN_+$:
\begin{equation}\label{piede}
\begin{split}
& \bigl|\int _0 ^M \frac{x^{\a-1} }{\l-x} ds \bigr| \geq c\,|\l|^{-2}, \qquad
 \forall \l\in \G_M \cup \g_\infty \\
&  \int _0 ^M \frac{x^{\a-1} }{|\l-x|} dx \leq c',\qquad \forall 
 \l\in \G_M \cup \g_\infty \\
& \int _M ^\infty \frac{x^{\a-1} }{|\l-x|} dx \leq c' M^{\a-1} \ln M ,\qquad \forall \l \in \g_\infty.
\end{split}
\end{equation}

From the  above estimates we infer
\begin{equation}\label{musica}
\bigl|
 g_M(t,t_w)- g(t,t_w) 
\bigr| \leq c\, M^{\a-1}\ln M
\end{equation}
where 
\begin{equation}\label{schianto}
g(t,t_w):= \frac{1}{2\pi i}\int_{\g_\infty } \frac{e^{- t_w \l }}{\l} 
\frac{ \int _0 ^\infty  \frac{e^{-x t }}{\l-x} x^{\a-1} dx }{
  \int _0 ^\infty \frac{1}{\l-x} x^{\a-1} dx } d\l
\end{equation}
Using the analytic properties of the integrand in the r.h.s. of
(\ref{schianto}),
one can show that $g(t,t_w)= g(t/t_w, 1)$. In order to compute $g(\theta, 1)$, we observe that for a suitable positive constant $c$
   $| g(\theta s,s)- g_M(\theta s , s)|\leq c M^{\a-1}\ln M$, for any $s\geq 1$ (in fact, the constant $c$ in (\ref{musica}) can be chosen uniformly if $t_w\geq 1$). 
By the results of Subsection \ref{step1} (compare (\ref{marcopolo}) with 
$\Pi(t,t_w) $ in Proposition \ref{schneider}) 
we get 
$$
\lim _{s\uto\infty} g_M (\theta s , s) = \text{ r.h.s. of } 
(\ref{noiosetto})
$$
thus concluding the proof.
\qed

%
%

\section{Other correlation function when $\t_0\dto 0$ after $E\dto 
-\infty$ }\label{pantani}

As stated in Proposition \ref{mamma}, the standard time--time correlation function $\Pi_E (t,t_w)$ has trivial behaviour after taking the limits $E\dto -\infty$, $\t_0\dto 0$. For physical reasons, it is more natural to disregard jumps into states with physical waiting time
$T_i=e^{E_i}$  much smaller than $\t_0$, since we assume that the physical instruments in the laboratory have sensibility up to the time unit $\t_0$. 
Therefore, let us fix $\d>0$ and  consider here  the more natural 
time--time  correlation function
$$
 \P^{(1)}_E(t,t_w)=  
\PP_E \bigl(\, x_E(u) \geq \d  \quad \forall u\in (t_w,t_w+t]:\;
 x_E(u)\not = x_E(u^-)\;\bigr).
$$
The main result of this section if the following one: 
\begin{prop}\label{marcopantani1000}
For almost all $\underline E$ 
\begin{equation}\label{corr-equiv2}
\lim _{t_w \uto\infty}\sup_{t>0} 
\varlimsup _{E_0\dto -\infty } \varlimsup _{E\dto -\infty}
 \bigl| \,  \P^{(1)}_E(  t ,t_w)-  \P_E(t ,t_w )\,\bigr|=0.
\end{equation}
In particular, for almost all $\underline E$,
\begin{equation}\label{terranera}
 \lim _{t_w\uto\infty}
 \lim_{E_0\dto -\infty}
\lim_{E\dto - \infty} \Pi_E^{(1)}\bigl(\theta  t_w,t_w\bigr)=
  \frac{\sin (\pi \a)}{\pi}  \int_{\frac \theta {1+\theta }}^1 u^{-\a}(1-u)^{\a-1}du, \qquad \forall \theta>0.
\end{equation}
\end{prop}

Note that the correlation function $\Pi^{(1)}_E(t,t_w)$ is the analogous of
$\Pi^{(1)}_N(t,t_w)$ of Proposition \ref{rahel}.  As for the proof of Proposition \ref{rahel}, a useful observation is that, given $\d>0$, 
\begin{equation}\label{joseph2}
\lim_{t \uto\infty} \lim _{\t_0 \dto 0}\lim _{E\dto -\infty}
\PP _E\bigl(\, x_E(t)> \d \,\bigr)  = 0\qquad \text{a.s.} .
\end{equation}
We can prove a stronger result concerning  to the phenomenon that 
with high probability the system visits deeper and deeper traps. In fact, note that by
(\ref{caporale})
$$
\PP _E\bigl(\, x_E(t)>\d \,\bigr) =
\frac{1}{2\pi i}\int_{\g_E} \frac{ e^{-t\l}}{\l}
\frac{ \sum _{j=1} ^{N_E}  \frac{\II_{x_j\geq \d} }{x_j-\l} }{\sum_{j=1}^{N_E}  \frac{1}{x_j-\l}}d\l.
$$
Then, by reasoning as in the proof of Proposition \ref{annibale} and using the
results of Section \ref{step2}, one can easily show the analogous of
Proposition \ref{bruna}:
\begin{prop}\label{volareohoh}
For almost all energy landscapes $\underline E$,
\begin{equation}\label{polo2}
\lim_{t \uto\infty} \lim _{\t_0 \dto 0}\lim _{E\dto -\infty}
t^{1-\a} \PP _E\bigl(\, x_E(t)>\d\,\bigr)  = \frac{ B(\d) }{c(\a)}
\end{equation}
where 
$$
B(\d):=\frac{  \int _{\d} ^\infty x^{\a-2} dx}{ \int_0^\infty
\frac{x^{\a-1}}{1+x} dx },\qquad c(\a):= 
\int_0 ^\infty y^{\a-1}e^{-y} dy
$$
\end{prop}

\noindent
{\bf Proof of Proposition \ref{marcopantani1000}}. Trivially, (\ref{terranera}) follows  from (\ref{corr-equiv2}) and Proposition \ref{mamma}. 
 Since $E_0\dto-\infty$ after $E\dto-\infty$, we assume  that 
 $E<E_0$ and  define 
$$
D_{E,E_0}:=\{i\,: E\leq  E_i \leq E_0  \}.
$$
This set corresponds to the small traps where we allow the particle to jump 
in.
By (\ref{joseph2})  and the same arguments used in the proof of
Proposition \ref{rahel} for $i=1$ (with exclusion of the last step  since here 
$ \lim _{E\dto-\infty} \frac{|D_{E,E_0}|}{N_E} =1$ a.s.)
 it is simple to derive  the assertion of the proposition 
   from Lemma \ref{smetana}.
\noindent
\qed

\begin{lemma}\label{smetana} For almost all energy landscapes $\underline E$ 
there exist positive constants  $p, c$  (independent of $E,E_0$) 
satisfying the   following property.  Whenever 
 $|\{  i\,:\, E_i >E_0 \}|>0 $, 
\begin{equation}\label{gegendiewand}
\limsup_{E\dto-\infty } \frac{1}{N_E}
 \sum _{i\in D_{E,E_0}}  \varphi _{E,E_0}(i,t)\leq  c\, e^{-p\, t},
\qquad \forall t>0,
\end{equation}
where $D_{E,E_0}:=\{i\,: E\leq  E_i \leq E_0  \}$, for $E<E_0$, and 
the function $\ \varphi _{E,E_0}$ is defined as
\begin{equation}\label{elisir11}
 \varphi _{E,E_0}(i,u):= 
\PP_E (Y_E(u)\in D_{E,E_0}\;\forall u\in [0,s]\,|\, Y_E(0)=i),
\end{equation}
\end{lemma}
\begin{proof} Let us assume that  $E<E_0$,  $|\{  i\,:\, E_i >E_0 \}|>0 $
and, without loss of generality, $\d=1$.

We  fix $\ell>0 $ such that $ e^{-\a\ell}<\frac{1}{2}$ and define 
\begin{align*}
& W_{1,E} :=\{i\,:\, E\leq E_i <E+\ell \},\; \;N_{1,E}:= |W_{1,E}|\\
&  W_{2,E} :=\{i\,:\, E+\ell \leq E_i \leq E_0 \},\;\;
 N_{2,E}:=|W_{2,E} |
\end{align*}
Note that $D_{E,E_0}=W_{1,E}\cup  W_{2,E}$ and that 
  $N_E, N_{1,E}, N_{2,E}$ are  Poisson variables having  expectations
$ e^{-\a E}$, $e^{-\a E}(1-e^{-\a\ell})$ and $e^{-\a E-\a\ell} -e^{-\a E_0}$.
In particular,  for almost all $\underline E$,
\begin{equation}\label{carofabri}
\lim _{E\dto-\infty} \frac{N_E }{e^{-\a E} } = 1,\qquad
p_1:=\lim _{E\dto-\infty} \frac{N_{1,E} }{N_E} = 1-e^{-\a\ell},\qquad
p_2:=\lim _{E\dto-\infty} \frac{N_{2,E} }{N_E} = e^{-\a\ell}<\frac{1}{2}
\end{equation}
We observe  that $ \tilde n:=N_E-N_{1,E} -N_{2,E}$ is a positive integer
 independent of $E$ and $x_i\geq e^{E_0-E-\ell}$ if $i\in W_{1,E}$, while 
$x_i\geq 1$ if $i\in W_{1,E}$. Let us introduce a new random walk 
$Y^\ast_E(t)$ on $\cS_E$ whose infinitesimal generator $\LL_E ^\ast$
is defined as $\LL_E$ with $x_i$ replaced by $x ^\ast _i$ defined as 
$$
x^\ast _i= 
\begin{cases} 
 x_i & \text{ if } i\not \in D_{E,E_0} ,\\
 A_E:=e^{E_0-E-\ell} &\text{ if } i\in W_{1,E} ,\\
1 &\text{ if } i\in W_{2,E}
\end{cases}
$$
We denote by   $\PP^\ast_E$ the probability on path space associated to $Y^\ast_E(t)$ when having initial uniform distribution and we set
\begin{equation}\label{elisir12}
 \varphi ^\ast _{E,E_0}(i,u):= 
\PP^\ast _E (Y^\ast _E(u)\in D_{E,E_0}\;\forall u\in [0,s]\,|\, 
Y ^\ast _E(0)=i).
\end{equation}
By a simple coupling argument, one gets
$\varphi _{E,E_0}(i,u)\leq  \varphi ^\ast _{E,E_0}(i,u)$. In particular
\begin{equation}\label{abbruzzo}
\frac{1}{N_E}\sum _{i\in D_{E,E_0} }
 \varphi _{E,E_0}(i,t)\leq \Phi:=\frac{1}{N_E}\sum _{i\in D_{E,E_0}}
 \varphi ^\ast _{E,E_0}(i,t),\qquad 
\forall i\in D_{E,E_0}, \forall t\geq 0.
\end{equation}

At this point it remains to estimate $\Phi$. In order to simplify notation we 
write simply $D,N,N_1,N_2,A$, by dropping the index $E$. Moreover, we consider the following realization of the dynamics of $Y^\ast_E$: after arriving at a site $i$, the system waits an exponential time of parameter $x_i^\ast$ and the it jumps to a point of $\cS_E$ with uniform probability. In particular, jumps can be \emph{degenerate}, i.e.  initial and final sites can coincide.

We claim that 
$$
 \Phi = \Phi _1+\Phi_2 +\Phi_3
$$
where 
\begin{align}
& \Phi_1:= \sum_{k_1=1}^\infty \sum_{k_2=0}^\infty
 {k_1+k_2 \choose k_1}
\bigl(\frac{N_1}{N}\bigr)^{k_1}\bigl(\frac{N_2 }{N}\bigr)^{k_2+1}
A^{k_1}\int _0^t d u\,  
e^{-Au}\frac{ u^{k_1-1}}{(k_1-1) !}
e^{-(t-u)}\frac{(t-u)^{k_2}}{k_2!}
\label{delphi1}\\
&\Phi_2 :=
\sum _{k_1=0}^\infty \sum_{k_2=1}^\infty
 {k_1+k_2 \choose k_1}
\bigl(\frac{N_1}{N}\bigr)^{k_1+1}\bigl(\frac{N_2 }{N}\bigr)^{k_2}
A^{k_1}\int _0^t d u\,  
e^{-Au}\frac{ u^{k_1}}{k_1 !}e^{-(t-u)}\frac{(t-u)^{k_2-1}}{(k_2-1)!} \label{delphi2}\\
&\Phi_3:=\frac{N_1}{N} \exp \bigl\{ -At\bigl(1-\frac{N_1}{N}
\bigr)\bigr\}+ \frac{N_2}{N} \exp \bigl\{
-t\bigl(1-\frac{N_2}{N}\bigr) \bigr\}
\end{align}
The above identities can be derived from the probabilistic interpretation of
 $k_1$, $k_2$ as 
$$ 
k_i=|\{\text{ jumps performed before time $t$ having starting point in } W_i \}|
$$
and from the following simple identities:
\begin{align*} 
& \PP\bigl(\, T_1+T_2+\cdots +T_n\in [z,z+dz)\,\bigr)=e^{-\k\,z} \k^n
\frac{ z^{n-1}}{(n-1)!}dz \\
&  \PP\bigl(\, T_1+T_2+\cdots +T_n \leq z\text{ and }
 T_1+T_2+\cdots +T_n+T_{n+1}>z  \,\bigr)=e^{-\k\,z} \k^n
\frac{ z^n}{n!}
\end{align*} 
where $z\geq 0$ and
 $T_1, T_2, \dots$ are  independent exponential variables with parameter $\k$.

Finally, we only need to prove  that,
 for suitable positive constant $c,p>0$,
\begin{equation}\label{PANTANI}
\limsup _{E\dto -\infty} \Phi_i \leq c e^{-pt},\qquad \forall i=1,2,3.
\end{equation}
We give the proof in the case $i=1$, the case $i=2$ is completely similar while the case $i=3$ follows directly from (\ref{carofabri}). \\

We fix $\g:\a<\g<1$,  set $k_0:= A^\g$ and  write
$$
\Phi_1= \Phi_1^{\leq k_0}+\Phi_1^{> k_0}
$$
where $ \Phi_1^{\leq k_0} $ is the contribution to $\Phi_1$ of addenda
in the r.h.s. of (\ref{delphi1}) with $1\leq k_1 \leq k_0$ and
 $k_2\geq 0$. \\

If $k_1> k_0$ then
$$
\Bigl( \frac{N_1}{N} \Bigr) ^{k_1} =
\Bigl( \frac{N-N_2}{N} \Bigr)^{k_1} 
\Bigl( 1-\frac{\tilde n}{N-N_2} \Bigr)^{k_1}\leq 
\Bigl( \frac{N-N_2}{N} \Bigr)^{k_1}  
 \Bigl( 1-\frac{1}{N}\Bigr)
^{A^\g}
$$
thus implying that
$$
\Phi_1^{> k_0}\leq \Phi_1 
\Bigl(
\frac{N-N_2}{N} ,N_2  
\Bigr) \Bigl( 1-\frac{1}{N}\Bigr)
^{A^\g} \leq  \Bigl( 1-\frac{1}{N}\Bigr)
^{A^\g} \dto 0 \qquad \text{ as } E\dto-\infty.
$$
where $\Phi_1\Bigl(\frac{N-N_2}{N} ,N_2  \Bigr)$
 is defined as in the r.h.s. of (\ref{delphi1}) with  $N_1$ replaced 
by $N-N_1$. Note that is  does not 
exceed $1$
 since it corresponds to the probability of a certain event.

Let us now consider the term  $\Phi_1^{\leq k_0}$. To this aim, since
$ {k_1+k_2 \choose k_1}\leq 2^{k_1 +k_2}$, 
$$
 \Phi_1^{\leq  k_0}\leq  \sum _{k_1=1}^{k_0}
 \sum_{k_2=0}^\infty
\bigl(\frac{2 A N_1}{N}\bigr)^{k_1}\bigl(\frac{2 t N_2 }{N}\bigr)^{k_2}
 \frac{I(0,t)}{(k_1-1)! k_2!}
$$
where
$$
I(w_1,w_2):=\int _{w_1}^{w_2} \,  e^{-(t-u) -Au}  u^{k_1-1} du.
$$
Fix $m>0$ with $ 2p_2+2mp_1<1$ (recall that $2p_2<1$). Then trivially
\begin{equation}\label{rimini}
I(0,\frac{m\,t}{A})\leq 
\frac{1}{k_1}  e^{-t\bigl(1-\frac{m}{A} \bigr) } 
\bigl(\frac{mt}{A} \bigr)^{k_1}  
\end{equation} 
From such a  bound, one gets immediately
\begin{equation}\label{adria}
\sum _{k_1=1}^{k_0}
 \sum_{k_2=0}^\infty
\bigl(\frac{2 A N_1}{N}\bigr)^{k_1}\bigl(\frac{2 t N_2 }{N}\bigr)^{k_2}
 \frac{1}{(k_1-1) ! k_2!} I( 0,\frac{m\,t}{A}) \leq
c\, \exp\Bigl\{ -t\Bigl(
 1-\frac{m}{A}-2\frac{N_1}{N} m -2\frac{N_2}{N} \Bigr)\Bigr\}
\end{equation}
The last expression, when $E\dto-\infty$, converges to
$c\exp\bigl( -t(1-2p_1 m -2p_2) \bigr)$, in agreement with
(\ref{PANTANI}).

\vspace{0.2cm}

In order to estimate the integral $I(\frac{m\,t}{A},t)=
  e^{-t} 
 \int _{\frac{m\,t}{A}  }^{t} \,  e^{- (A-1)u}  u^{k_1-1}$,
 we observe that  
$$
 \int_s ^w e^{-z u }u^n du = (-1)^n \frac{d}{dz^n}\frac{ e^{-zs}-e^{-zw} }{z}
\leq \frac{e^{-z s}}{z}  (s +\frac{1}{z} )^n+ \frac{e^{-z w}}{z}  (w +\frac{1}{z} )^n,\qquad \forall s,w,z\geq 0,
$$
thus implying the bound
$$
I (\frac{m\,t}{A},t)\leq c  e^{-t} A^{-k_1} (mt+1 )^{k_1 -1} +
e^{-At} (A-1)^{-1} (t+ \frac{1}{A-1}) ^{k_1-1}
$$
The contribution of  $ c  e^{-t} A^{-k_1} (mt+1 )^{k_1 -1}$ to 
$\Phi_1^{>k_0}$ can be treated by means of  estimates similar to the ones
 leading to  (\ref{adria}).

In order to conclude we only need to show that
\begin{equation}\label{tarzan}
e^{-At} \sum _{k_1=1}^{A^\g}
 \sum_{k_2=1}^\infty
\bigl(\frac{2 A t N_1}{N}\bigr)^{k_1-1}
\bigl(\frac{2 t N_2 }{N}\bigr)^{k_2}
 \frac{1}{(k_1-1)! k_2!}  \dto 0,\quad  \text{ as } E\dto-\infty.
\end{equation}
To this aim observe that
\begin{equation*}
\begin{split}
\text{r.h.s. of } (\ref{tarzan})
 & \leq 
e^{-At+ 2t  \frac{ N_2 }{N}} 
 \sum _{k_1=0}^{A^\g}
\bigl(\frac{2 A t N_1}{N}\bigr)^{k_1}\\
  &\leq c(t) e^{-At}
A^\g (4tA)^{A^\g}=c(t) \exp\Bigl\{-At +\g \ln A + A^\g \ln( 2t A) \Bigr\}
\end{split}
\end{equation*}
Since $0<\g<1$, we get (\ref{tarzan}), thus concluding the proof.
\end{proof}

\appendix

\section{Laplace Transform}

\begin{prop}\label{LTI}
Let  $G(t)$ be a bounded measurable function on
$(0,\infty)$ and let us consider the Laplace transform
$$
\hat G (\o) =\int _0 ^\infty G(t) e^{-t\o} dt
$$
well defined if $\Re(\o)>0$. Let us define
\begin{equation}\label{tigre}
\cA:= \{ re^{i\theta}\,:\, 0<r<\infty,\; |\theta|\leq \frac{3}{4}\pi \}
\end{equation}
Suppose that $\hat G$ can be analytically continued to 
$\CC \setminus (-\infty, 0]$ and that  there are   positive constants 
 $\g,\b,\a, c$ and   $B\in \RR$ such that 
\begin{align}
& |\,\hat G(\o)\,| \leq c\,|\o|^{-\g} \qquad \qquad \forall \o\in \cA,\; |\o|\geq 1,
\label{nico}\\
& \,|\o^\b \hat G(\o) -B\,|\leq c\,|\o|^{\a} \quad   \;
\forall \o\in \cA,\; |\o|\leq 1 \label{france} .
\end{align}
Then,
$$ \lim _{s\uto\infty} s^{1-\b} G(s) =\frac{B}{c(\b)}
\text{ where } c(\b) := \int_0 ^\infty y^{\b-1} e^{-y} dy .
$$
\end{prop}
\begin{proof}
 If we set $H(s):=\frac{B}{c(\b)} s^{\b-1}$ with    $s>0$, 
then the Laplace transform $\hat H(\o)$ is well defined for 
$\Re(\o)>0$, $\hat H(\o) = B \o^{-\b} $ and trivially   $\hat H(\o)$ 
can be analytically continued to $\CC \setminus (-\infty, 0]$.

By the inverse formula for Laplace transform (see Chapter 4, Section 4 in
\cite{doe1}), we have
\begin{equation}\label{dostoevski}
  G(s)= \lim_{K\to \infty}
  \frac{1}{2\pi i} \int _{x-iK} ^{x+iK}
  e^{s \o} \hat G(\o) d\o,
  \qquad \forall s>0,\; x>0,
\end{equation}
where $\o$ runs over the vertical path connecting $x-iK$ and $x+iK$. 
The above formula remains true if substituting $G$ with $H$. Therefore,
\begin{equation}\label{patroclo}
s^{1-\b} G(s) -\frac{B}{c(\b)} =
 \lim_{K\to \infty}
 \frac{ s^{1-\b} }{2\pi i} \int _{x-iK} ^{x+iK}
  e^{s \o}\bigl( \hat G(\o)-\hat H(\o) \bigr) d\o,
  \quad \forall s>0,\; x>0,
\end{equation}

Let $\rho:= \min(\g,\b)/2$. Fix a positive number $x$ and, 
given $K$ and $s$, define the following paths (see figure below).

$\g_K$ is the vertical path from  $x-iK$ to  $x+iK$.
$\g_{1,+}$ is the segment from $-s^{-1}+s^{-1}i$ to $-1+i$.
 $\g_{2,+}$ is an arc from    $-1+i$ to $-K^\rho +iK$ given 
by the parametrisation $z(t)=-t+it^{1/\rho}$ with $t\in [1,K^\rho ]$.
  $\g_{3,+}$ is the horizontal segment from  $-K^\rho +iK$ to
$x+iK$.
For $i=1,2,3$ we define the path $\g_{i,-}$ by considering 
the reflection of  $\g_{i,+}$ w.r.t. the real axis and inverting 
the orientation. Let  $\g_0$  be the positive--oriented circular
  arc of radius $s^{-1}$ from $-s^{-1}-s^{-1}i$ to
 $-s^{-1}+s^{-1}i$ crossing the axis of positive real numbers.\\
Note that the above paths depend on $s$ and/or $K$.

\begin{figure}[!ht]
    \begin{center}
       \psfrag{a}[l][l]{$x+iK$}
       \psfrag{b}[l][l]{$x-iK$}
       \psfrag{c}[l][l]{$-K^\rho+iK$}
        \psfrag{d}[l][l]{$-K^\rho-iK$}
       \psfrag{1}[l][l]{$\gamma_K$}
       \psfrag{2}[l][l]{$\gamma_{3,+}$}
       \psfrag{3}[l][l]{$\gamma_{2,+}$}
       \psfrag{4}[l][l]{$\gamma_{1,+}$}
       \psfrag{5}[l][l]{$\g_0$}
       \psfrag{7}[l][l]{$\g_{2,-}$}
       \psfrag{6}[l][l]{$\g_{1,-}$}
       \psfrag{8}[l][l]{$\hat \gamma_{3,-}$}
      \includegraphics[width=5cm]{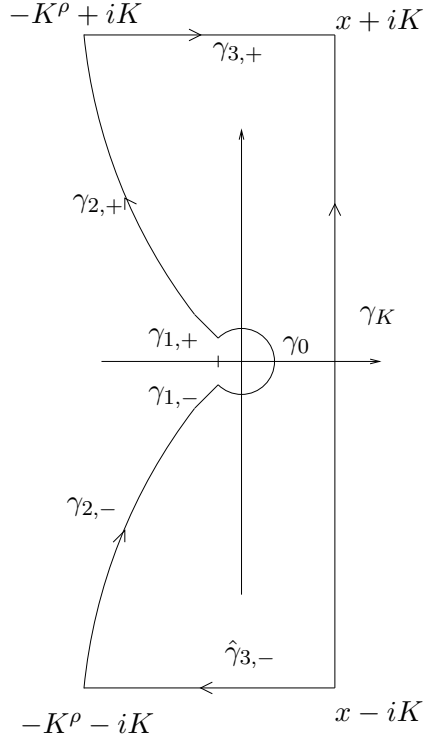}
      \caption{The integration paths}
    \label{f:1}
    \end{center}
  \end{figure}

Because of analyticity,  the integral over  $\g_K$ of   
$ e^{s \o} \hat G _(\o) $ is equal to the sum of the integrals 
over $\g_{3,-}$, $\g_{2,-}$, $\g_{1,-}$, $\g_0$,  $\g_{1,+}$, 
$\g_{2,+}$, $\g_{3,+}$. The same is valid with $\hat G$ 
replaced with  $\hat H$.\\

By  (\ref{nico}) we have that
\begin{equation}
\int_{\g_{3,\pm} } |e^{s \o} \hat G(\o) | \, |d\o| 
\leq c\, e^{sx} K^{\rho-\g}\dto 0\qquad \text{ as } K\uto \infty
\end{equation}
and, for a suitable rational function $f$,
\begin{equation}\label{help}
\begin{split}
 s^{1-\b}\int_{\g_{2,\pm} } |e^{s \o} \hat G(\o) | \, |d\o| & \leq
s^{1-\b} \int _1 ^\infty \bigl|e^{s(-t+it^{1/\rho})} \hat 
 G( -t+it^{1/\rho} )
(-1+\rho^{-1}t^{  \frac{1-\rho}{\rho}} i)\,\bigr |dt\\
& \leq  c\,s^{1-\b} \int_1^\infty e^{-st} f(t) dt\leq c' s^{1-\b} 
e^{-\frac{s}{2}}\dto 0 , \qquad\text{ as } s\uto \infty .
\end{split}
\end{equation}
Similarly, it can be proved that the corresponding integrals with 
$\hat G$ substituted with $\hat H$ go to $0$ by taking the limits 
$K\uto \infty$, $s\uto\infty$.

Let us now estimate $s^{1-\b} \int_{s^{-1}} ^1 e^{-st} t^{\a-\b} dt$
by dividing the path of integration in two paths. Choosing $0<\d<1$,
\begin{equation*}
\begin{split}
& s^{1-\b} \int_{s^{-1}} ^{ s^{-1+\d} }  e^{-st} t^{\a-\b} dt \leq  c\,
s^{1-\b} \bigl| s^{(-1+\d)(1+\a-\b)} -  s^{-(1+\a-\b)}\bigr|
\leq c' s^{-\a} +c ' s^{-\a+\d(1+\a-\b)},   \\
& s^{1-\b} \int_{s^{-1+\d}} ^1   e^{-st} t^{\a-\b} dt \leq
s^{1-\b}  e^{-s^\d} g(s)
\end{split}
\end{equation*}
for a suitable rational function $g(s)$. In particular, choosing $\d$  
small enough, the above upper bounds imply
$$
\lim _{s\uto\infty} s^{1-\b} \int_{s^{-1}} ^1 e^{-st} t^{\a-\b} dt=0.
$$
This result, together with (\ref{france}), implies
$$
\lim _{s\uto\infty}
 s^{1-\b}  \frac{1}{2\pi i} \int _{\g_{1,\pm} }
  e^{s \o}\bigl( \hat G(\o)-\hat H(\o) \bigr) d\o =0.
$$
Trivially, by (\ref{france}),
$$
 s^{1-\b}  \bigl | \int _{\g_0 }
  e^{s \o}\bigl( \hat G(\o)-\hat H(\o) \bigr) d\o\bigr|\leq
s^{-\a}\dto 0 \qquad\text{ as } s\uto \infty .
$$
The Proposition now  follows from the   estimates above and (\ref{patroclo}).
\end{proof}

\section{Perturbation Theory}\label{step5}

In this appendix  we comment  
on a paper by Melin and Butaud \cite{BM} where the eigenvalues and
eigenfunctions of the generator of our model were computed using
perturbation theory. As pointed out earlier, these results are at
variance with our exact results, and it may be worthwhile to point out
the flaw in their arguments.
Melin and Butaud  write the
generator $\LL$ defined in (\ref{perleaiporci}) as
$\LL=T+\frac{1}{N}  T^{(1)} $ where
\begin{equation}
T :=
\begin{pmatrix}
x_1 & 0 &\ldots & 0 \\
0 & x_2 & \ldots & 0 \\
\vdots        & \vdots        & \ddots & \vdots \\
0  & 0  & \ldots &  x_N
\end{pmatrix},\qquad
 T^{(1)}:=
\begin{pmatrix}
-x_1 & -x_1 &\ldots & -x_1 \\
-x_2 & -x_2 & \ldots & -x_2 \\
\vdots        & \vdots        & \ddots & \vdots \\
-x_N  & -x_N  & \ldots & -x_N
\end{pmatrix}
\end{equation}
The factor $1/N$ in front of the second term encourages them to 
consider this term as a small perturbation. 
 Both $T$ and $T^{(1)}$ are symmetric operators on  $L^2(\mu)$ 
where $\mu(i):= x_i^{-1}$. We denote $<\cdot,\cdot>$  
the scalar product in   $L^2(\mu)$ and  assume $x_1, \dots , x_N$ to be
distinct positive numbers.

Given an operator $A: L^2(\mu)\to L^2(,\mu)$ we write $\|A\|$ for
its operator norm.
 Because of symmetry, $\|T\|$ and $\|T^{(1)}\|$ are given 
by the maximum of $|\l|$, with $\l$ eigenvalue. Trivially, 
$T$ has eigenvalues $x_1,\dots , x_N$ and $T(e_i)=x_i e_i$ 
where $e_1, \dots ,e_N$ is the
canonical basis of $\RR^N$, while
 $T^{(1)}$ has eigenvalues $0,-(x_1+x_2+\dots +x_N)$.

Given $z\in \CC$ we can define the holomorphic function $ T(z)=
T+zT^{(1)}$. A natural condition in order to apply perturbation
theory to $T(z)$  (see \cite{Kato}, chapter II)  is
\begin{equation}\label{tatatata}
 |z|< \frac{d}{2 a_0}
\end{equation}
where
$$
 d=\inf_{i\not = j}|x_i-x_j|, \qquad a_0 := \min _{a\in \RR}\|T^{(1)}-a\|=
\frac{x_1+x_2+\dots +x_N}{2}.
$$
In this case, we can conclude that $ T(z)= T+zT^{(1)}$ has $N$
eigenvalues $\l_1(z),\dots,\l_N(z)$ with
$\l_k(z)=\sum_{n=0}^\infty \l^{(n)}_k z^n$,  $|\l_k ^{(n)} | \leq
a_0 ^n (2/d)^{n-1} $  and
\begin{align*}
&\l^{(1)}_k=\frac{<T^{(1)} e_k, e_k >}{<e_k,e_k>}\\
&\l^{(2)}_k=\sum_{j\not=k} (x_k-x_j)^{-1}
\frac{<T^{(1)} e_k, e_j >^2}{<e_k,e_k><e_j,e_j>}\\
& ...
\end{align*}
Similar series exist for  the perturbed eigenvectors.

However,  the crucial condition (\ref{tatatata}) is hardly  
satisfied when  $z=\frac{1}{N}$, since it reads
\begin{equation}\label{quelle}
 \Av _{j=1}^N x_j \leq \inf _{ i\neq j }|x_i-x_j|
\end{equation}
while a.s.  the l.h.s. of (\ref{quelle}) has non zero limit and
the r.h.s. converges to $0$ like $1/N$. 

The fact that the conditions for the application of perturbation
theory are violated explains why its predictions are incorrect. 
This discrepancy happens not to be too obvious as far as the
eigenvalues are concerned (which are caught between the diagonal
elements of the generator and thus are somewhat similar to them,
but the shape of the eigenfunctions is sharply different).

Namely, by of Proposition $\ref{structure1}$,  when $j\not =1$ and  $
x_{j-1},x_j$ are very near each other, the eigenvector
$\psi^{(j)}$ related to the eigenvalue $\l_j: x_{j-1}<\l_j<x_j$
has two main peaks of opposite sign given by $\psi^{(j)}_{j-1}$
and $\psi^{(j)}_j$, this is very different from the predictions  of \cite{BM}
(see their Figure 4).

\section{Complex integral representation}\label{rinnovo}
Let $\LL$ be a Markov generator on the state space  $\cS:=\{1,2,\dots, N \}$, 
reversible w.r.t. a positive measure $\mu$. We can think of $\LL$ as a linear operator
on $\RR^N$, symmetric w.r.t. to the scalar product $(\cdot,\cdot)_\mu$ where
$$ 
(a,b)_\mu=\sum _{i=1}^N \mu(i) a_i b_i 
$$
In what follows we endow  $\RR^N$  with the scalar product $(\cdot,\cdot)_\mu$
 (and not with the standard Euclidean scalar product).
Since $\LL$ is symmetric, we can  orthogonally decompose  $\RR^N$ as
$
\RR^N=W_1\oplus W_2\oplus \cdots \oplus W_m
$ 
such that $\LL=\sum _{k=1}^m \l_k P_{W_k}$, where $P_{W_k}$ denotes the orthogonal projection  on $W_k$ and $\l_i\not = \l_j$ if $i\not = j$. 
Given  $\l\in \CC \setminus \{\l_1,\dots,\l_m\}$, we write $R(\l)$ for the resolvent 
$$
R(\l):=(\l\II-\LL)^{-1}=\sum _{k=1}^m \frac{1}{\l-\l_k}  P_{W_k} .
$$ 
Then, the Residue Theorem implies the 
integral representation
\begin{equation}\label{ven11}
 e^{-t\LL}=\frac{1}{2\pi i}\int_\g  e^{-t\l} R(\l) d\l
\end{equation}
where $\g$ is a positive oriented loop containing in its interior $\l_1,\l_2,\dots,\l_m$.

Given a probability measure  $\nu$ on $\cS$,
we denote by $\PP_\nu$ the probability measure on the 
path space associated to the
continuous--time random walk $Y(t)$ on $\cS$ with generator $\LL$ and initial distribution $\nu$. Fix $j\in \cS$ and let  $v \in\RR^N $ be 
 such that $v_i=\d_{i,j}$. We write $\frac{d \nu}{d\mu}$ for the Radon derivate, i.e.  $\frac{d \nu}{d\mu} (i)= \frac{\nu(i)}{\mu(i)}$. Then the symmetry of $\LL$ w.r.t. the scalar product $(\cdot,\cdot)_\mu$ implies
\begin{equation}\label{ven12}
\begin{split}
\PP_\nu(\, Y(t) = j\,)& =\sum_{k=1}^N \nu(k) \bigl( e^{-t\LL } \bigr) _{k,j} =
\mu\bigl ( \frac{d \nu}{d\mu}, e^{-t\LL} v  \bigr) =
\mu\bigl ( e^{-t\LL} \frac{d \nu}{d\mu}, v \bigr)\\
& = \mu(j) \sum _{k=1}^N
 \bigl( e^{-t\LL} \bigr)_{j,k} \frac{\nu (k)}{\mu (k)}.
\end{split}
\end{equation}
By plugging (\ref{ven11}) in the r.h.s. of  (\ref{ven12})
 we get the  integral representation
\begin{equation}\label{lupo}
\PP_\nu(\, Y(t) = j\,)=\frac{1}{2\pi i}\int_\g  e^{-t\l} \bigl\{
\sum _{k=1} ^N  \mu(j)  R_{jk}(\l)  \frac{\nu(k)}{\mu(k)} \bigr\}d\l.
\end{equation}
In particular, given $h$ function on $\cS$
$$
\EE_{\PP_\nu} \bigl(\,h(Y(t)\,\bigr) = \sum _{j=1}^N \sum_{k=1}^N
 \frac{1}{2\pi i}
h(j)  \mu(j) \frac{\nu(k)}{\mu(k)}  \int_\g  e^{-t\l} R_{jk}(\l) d\l
$$
thus allowing to get an integral representation for
$ 
\Pi (t,t_w):= \PP_\nu(\text{ no jump in } [t_w,t_w+t])
$.
If we set $\mu(i)=\t_i =x_i^{-1}$ and $\nu(i)=N^{-1}$ (uniform initial probability), then
\begin{equation}\label{jiribravo}
\PP_\nu(\, Y(t) = j\,)=\frac{1}{2\pi i}\int_\g  e^{-t\l} \bigl\{
\frac{1}{N} \sum _{k} \frac{x_k}{x_j}  R_{jk}(\l)  \bigr\}d\l.
\end{equation}\\

Let us consider now the  special case given by the Bouchaud's REM--like trap 
model
 where $\LL:=\LL_N$ is  defined in (\ref{perleaiporci}) and $\nu$ is the uniform distribution on $\cS$.
Note the all the  integral formulas obtained in Section \ref{inizio} can be 
derived from the following one:
\begin{equation}\label{babbobuono}
 \PP_\nu(\, Y(t) = j\,)=\frac{1}{2\pi i}\int_\g 
 e^{-t\l} \frac{1}{(\l-x_j)\phi(\l)}d\l
\end{equation}
where $\phi (\l)=\sum_{k=1}^N \frac{\l}{\l-x_k}$.
In what follows we prove that  (\ref{babbobuono}) corresponds to  (\ref{jiribravo}).

We know already that $\det(\l\II-\LL)$
 has distinct zeros given by the $N$  distinct zeros  of $\phi (\l)$. In particular, it must be
\begin{equation}\label{pig}
\det(\l \II -\LL ) =
\frac{1}{N} \phi (\l) \prod _j (\l-x_j)=
\frac{1}{N} \l \sum _k \prod  _{j:\,j\not = k} (\l-x_j) 
\end{equation}
Given a matrix $A$ we write $[A]_{i,j}$ for the determinant of the matrix obtain from $A$ by erasing the $i$--th row and the $j$-- column.
Since 
$$
 R_{j,k}(\l) = (-1)^{j+k+1} \frac{[\l \II-\LL ]_{k,j} }{\det(\l\II-\LL) }
$$
and due to (\ref{pig}), in order to derive  (\ref{babbobuono}) from 
(\ref{jiribravo}) we only have to show that
\begin{equation}\label{davideaiuta}
\sum _k (-1)^{j+k+1}\frac{ x_k }{x_j} [\l\II -\LL]_{k,j}=\prod _{s:\,
 s \not = j} (\l-x_s)
\end{equation}
In order to prove the above identity observe that  $[\l\II -\LL]_{k,j}$ is a polynomial of degree $N-1$ if $k=j$, otherwise it has degree $N-2$. The l.h.s. 
of (\ref{davideaiuta}) is a monomic polynomial of degree $N-1$. At this point we only have to verify that $x_s$,  $s\not = j$, are zeros of the l.h.s. 
of (\ref{davideaiuta}).
This is trivial if one observes that   the l.h.s. 
of (\ref{davideaiuta}) is the determinant of the matrix obtain from
$\l\II-\LL $ by replacing the $j$--column with the vector $w$ with 
$w_i=\frac{x_i}{x_j}$ for $i=1,2,\dots,N$. It is simple to verify that, 
if  $\l=x_s$ for some
$s\not =j $, the $j$--th row and the $s$--th row in such a matrix are proportional, thus implying the thesis. 

\vspace{1cm}

\noindent
{\bf Acknowledgements.} The authors thank  J. \v Cern\'y and V.  Gayrard
  for  useful discussions.


\end{document}